\documentclass[12pt]{article} 
\usepackage{latexsym, 
amscd, 
amsthm, 
graphicx, epsfig, amssymb}

\newtheorem{thm}{Theorem}[section]
\newtheorem{defn}[thm]{Definition}
\newtheorem{prop}[thm]{Proposition}
\newtheorem{cor}[thm]{Corollary}
\newtheorem{lemma}[thm]{Lemma}
\newtheorem{rema}[thm]{Remark}

\newcommand{\halmos}{\rule{1ex}{1.4ex}}

\newcommand{\nn}{\nonumber \\}

 \newcommand{\res}{\mbox{\rm Res}}

\renewcommand{\hom}{\mbox{\rm Hom}}

 \newcommand{\pf}{{\it Proof.}\hspace{2ex}}
 \newcommand{\epf}{\hspace*{\fill}\mbox{$\halmos$}}
 \newcommand{\epfv}{\hspace*{\fill}\mbox{$\halmos$}\vspace{1em}}

\newcommand{\wt}{\mbox{\rm wt}\ }

\newcommand{\lbar}{\bigg\vert}

\newcommand{\C}{\mathbb{C}}
\newcommand{\Z}{\mathbb{Z}}
\newcommand{\R}{\mathbb{R}}

\newcommand{\N}{\mathbb{N}}

\newcommand{\one}{\mathbf{1}}

\title{ {\bf Generalized twisted modules associated to general
automorphisms of a vertex operator algebra} }
\date{}
\author{Yi-Zhi Huang}

\begin{document}

\bibliographystyle{alpha}
\maketitle

\begin{abstract}
We introduce a notion of strongly $\C^{\times}$-graded, or equivalently,
$\C/\Z$-graded
generalized $g$-twisted $V$-module associated to an automorphism 
$g$, not necessarily of finite order, 
of a vertex operator algebra. We also introduce a notion of
strongly $\C$-graded
generalized $g$-twisted $V$-module 
if $V$ admits an additional $\C$-grading compatible with $g$.
Let $V=\coprod_{n\in \Z}V_{(n)}$ be a vertex operator 
algebra such that $V_{(0)}=\C\one$ and $V_{(n)}=0$ for
$n<0$ and let $u$ be an element of $V$ of weight $1$
such that  $L(1)u=0$.
Then the exponential of 
$2\pi \sqrt{-1}\; \res_{x} Y(u, x)$ is an 
automorphism $g_{u}$ of $V$. In this case, a strongly $\C$-graded
generalized $g_{u}$-twisted $V$-module 
is constructed from a strongly $\C$-graded generalized 
$V$-module with a compatible action of $g_{u}$ 
by modifying the vertex operator map for the generalized $V$-module 
using  the exponential of the negative-power 
part of the vertex operator $Y(u, x)$. In particular, we give 
examples of such generalized twisted modules associated to 
the exponentials of some screening operators on 
certain vertex operator algebras related to the triplet $W$-algebras. An 
important feature is that we have to work with generalized 
(twisted) $V$-modules which are doubly graded by the group
$\C/\Z$ or $\C$ and by generalized 
eigenspaces (not just eigenspaces) for $L(0)$, 
and the twisted vertex operators 
in general involve the logarithm of the formal variable. 
\end{abstract}

\renewcommand{\theequation}{\thesection.\arabic{equation}}
\renewcommand{\thethm}{\thesection.\arabic{thm}}
\setcounter{equation}{0}
\setcounter{thm}{0}

\section{Introduction}

The present paper is the first in a series of papers 
developing systematically a general theory of twisted modules
for vertex operator algebras. This 
theory is in fact equivalent to the study of orbifold 
theories in conformal field theory. Though twisted modules
associated to finite-order automorphisms of a vertex 
operator algebra have been introduced and studied (see below for 
more discussions and references), 
curiously, twisted modules
associated to infinite-order automorphisms of a vertex operator
algebras have not been previously formulated and studied mathematically.
To develop our general theory of twisted representations of 
vertex operator algebras, it is necessary to first 
find the correct category of twisted modules associated to both 
finite- and infinite-order automorphisms of a vertex operator algebra, and 
to study and construct such twisted modules. 
In this paper, we introduce a notion of strongly $\C^{\times}$-,
or equivalently,
$\C/\Z$-graded
generalized $g$-twisted module associated to an automorphism
of a vertex operator algebra
and more generally, a notion of 
strongly $\C$-graded
generalized $g$-twisted $V$-module when 
$V$ has an additional $\C$-grading compatible with $g$.
We also construct examples. The category of such
generalized twisted modules will be the category of twisted modules
that we shall study in this series of papers. 

Twisted modules for vertex operator algebras were introduced
in the construction of 
the moonshine module vertex operator algebra $V^{\natural}$
by I. Frenkel, J. Lepowsky and A. Meurman \cite{FLM1}, 
\cite{FLM2} and
\cite{FLM3}.  Given a vertex operator algebra $V$
and an automorphism $g$ of $V$ of finite order,  the
notion of $g$-twisted $V$-module is a natural but very subtle
generalization of the notion of $V$-module such that the action of 
$g$ is incorporated. The main axiom in the definition of 
the notion of $g$-twisted $V$-module is the twisted Jacobi
identity of the type proved in \cite{Le1}, \cite{FLM2}, \cite{Le2} 
and \cite{FLM3}. In fact, $V^{\natural}$
was the first example of an orbifold theory, before orbifold 
theories were studied in physics.

Orbifold theories are important examples of conformal 
field theories obtained from known conformal field theories and their 
automorphisms (see 
\cite{DHVW1}, \cite{DHVW2}, \cite{H}, \cite{DFMS}, \cite{HV}, \cite{NSV}, \cite{M},
\cite{DGH}, \cite{DVVV}, \cite{DGM}, as
well as, for example
\cite{KS}, \cite{FKS}, \cite{Ba1}, \cite{Ba2}, \cite{BHS},
\cite{dBHO}, \cite{HO}, \cite{GHHO}, \cite{Ba3} and \cite{HH}). 
Mathematically, the study of orbifold theories 
is equivalent to the theory of twisted 
representations of vertex operator algebras,
in the sense that any result or conjecture on orbifold conformal 
field theories can be formulated precisely as a result or conjecture
in the theory of twisted representations of vertex operator algebras.

There have been a number of papers constructing 
and studying twisted modules associated to finite-order automorphisms of
vertex operator algebras 
(see for example, \cite{Le1}, \cite{FLM2}, 
\cite{Le2}, \cite{FLM3}, \cite{D}, \cite{DL},
\cite{DonLM1}, \cite{DonLM2},  \cite{Li}, \cite{BDM}, \cite{DoyLM1},
\cite{DoyLM2}, \cite{BHL}, and the references in these papers, especially 
in the last three papers).
However, as far as the author knows,
the precise notion of $g$-twisted $V$-module for an infinite-order
automorphism $g$
has not been previously formulated, and such twisted modules have not 
been previously constructed and studied. 

The main difficulty in formulating such a general notion of $g$-twisted
$V$-module is that 
the right-hand side 
$$\frac{1}{k}\sum_{j\in \Z /k \Z}x_2^{-1}
\delta\left(\eta^j\frac{(x_1-x_0)^{1/k}}{x_2^{1/k}}\right)Y^g (Y(g^j
u,x_0)v,x_2)$$
of the twisted Jacobi identity (see \cite{Le1}, \cite{FLM2},
\cite{Le2} and \cite{FLM3};
cf. \cite{DL})
and the form of the expression
$$(x_{2}+x_{0})^{N}Y^{g}(Y(u, x_{0})v, x_{2})$$
needed in the 
weak associativity relation (see \cite{Li})
cannot be generalized to the general case.
In the present paper, we first reformulate the
notion of $g$-twisted $V$-module in the case of finite order $g$ using
a duality property (incorporating both commutativity and 
associativity) formulated in terms of complex variables
and then generalize this
equivalent definition to the general case. 
The precise generalization involves more subtle issues which 
fortunately have been treated in the logarithmic 
tensor product theory developed by Lepowsky, Zhang and the author
in \cite{HLZ1} and \cite{HLZ2}. First, in the general case, 
we have to consider a twisted generalization of the notion of generalized 
$V$-module in the sense of \cite{HLZ1} and \cite{HLZ2}
instead of the notion of $V$-module. 
Second, the twisted vertex operators 
in general involve the logarithm of the formal variable.
Third, we have to consider an additional grading by $\mathbb{C}^{\times}$
(equivalently by $\C/\Z$) or by $\C$
on such a twisted generalization to introduce a notion 
of strongly $\C/\Z$- or $\C$-graded generalized $g$-twisted 
$V$-module in the spirit of the 
the notion of strongly $\tilde{A}$-graded generalized $V$-module
in \cite{HLZ1} and \cite{HLZ2} for a suitable abelian group $\tilde{A}$.   
It is important to note
that in the general case, all these three new ingredients are 
necessary. In particular, many of the general results obtained in 
\cite{HLZ1} and \cite{HLZ2} are needed for 
the study of strongly $\mathbb{C}/\Z$- or $\C$-graded generalized $g$-twisted 
$V$-modules introduced and constructed in the present paper.

Even when  $g$ is a finite-order automorphism 
of $V$, how to explicitly construct a nontrivial $g$-twisted $V$-module
is still an open problem in general. The problem of constructing 
strongly
$\C/\Z$- or $\C$-graded generalized $g$-twisted
$V$-module for a general automorphism of $V$ is certainly more
difficult. However, when we consider certain special types of 
automorphisms of $V$, we might still be able to construct 
such generalized twisted
$V$-modules
associated to these automorphisms even though the general
problem is  not solved; indeed, in the finite-order case, the 
construction of certain twisted modules for lattice vertex operator
algebras was given by Frenkel, Lepowsky and Meurman
as one of the many hard steps in the construction 
of the moonshine module in the work \cite{FLM3}
and has been used to solve many 
other significant problems.

Given $u\in V_{(1)}$ such that  $L(1)u=0$, 
we show that
$g_{u}=e^{2\pi \sqrt{-1}\;
 \mbox{\rm \scriptsize $\res$}_{x} Y(u, x)}$ is an automorphism of 
$V$. In general, its order might  not be finite. 
In the present paper, generalizing a construction by Li 
\cite{Li} in the case that $\res_{x} Y(u, x)$
acts on $V$ semisimply and has eigenvalues belonging to $\frac{1}{k}\Z$,
we construct
a strongly
$\C$-graded generalized 
$g_{u}$-twisted $V$-module
for $u\in V_{(1)}$ from a strongly
$\C$-graded generalized $V$-module with an action
of $g_{u}$ such that 
the $\C$-grading is given by the generalized 
eigenspaces of the action of $\res_{x} Y(u, x)$.
We modify the vertex operator map of the generalized $V$-module
using the exponentials of the negative-power
part of the vertex operator $Y(u, x)$. 
We study such exponentials for $u$ not necessarily of weight $1$
and apply the result to the case of $u\in V_{(1)}$ later in the 
construction.
Our results show the importance of studying exponentials of 
vertex operators. 

As explicit examples, we apply this construction to 
certain vertex operator algebras constructed from 
certain lattice vertex operator algebras and 
automorphisms obtained by exponentiating
suitable screening operators. The fixed-point subalgebras of these vertex operator
algebras under the group generated by 
these automorphisms are in fact the triplet $\mathcal{W}$-algebras which in recent years
have attracted a lot of attention from physicists and mathematicians.
These algebras were
introduced first by Kausch \cite{K1} and have been studied extensively
by Flohr \cite{F1} \cite{F2}, 
Gaberdiel-Kausch \cite{GK1} \cite{GK2}, Kausch \cite{K2}, 
Fuchs-Hwang-Semikhatov-Tipunin \cite{FHST}, Abe \cite{A},
Feigin-Ga{\u\i}nutdinov-Semikhatov-Tipunin \cite{FGST1} \cite{FGST}
\cite{FGST3}, Carqueville-Flohr \cite{CF}, Flohr-Gaberdiel 
\cite{FG}, Fuchs \cite{Fu}, Adamovi\'{c}-Milas \cite{AM1} \cite{AM2} \cite{AM},
Flohr-Grabow-Koehn \cite{FGK}, 
Flohr-Knuth \cite{FK}, Gaberdiel-Runkel \cite{GR1} \cite{GR2}, 
Ga{\u\i}nutdinov-Tipunin \cite{GT}, Pearce-Rasmussen-Ruelle \cite{PRR1}
\cite{PRR2}, Nagatomo-Tsuchiya
\cite{NT} and Rasmussen \cite{R}.
Our results
show, in particular, that in this special case, the intertwining operators 
for suitable fixed-point subalgebras 
given by Adamovi\'{c} and Milas in Theorem 9.1 in \cite{AM} 
are in fact twisted vertex operators.
This connection
between suitable intertwining operators and twisted vertex operators 
shows that the triplet logarithmic conformal field theories 
are closely related to orbifold theories. 
Note that up to now, the triplet $W$-algebras are the main 
known examples of vertex operator algebras 
for logarithmic conformal field theories. We expect
that our theory will provide an orbifold-theoretic 
approach to the triplet logarithmic conformal field theories
and will give more interesting examples of logarithmic 
conformal field theories. 

There are certainly many examples of vertex operator 
algebras which have automorphisms of infinite orders. 
For example,  the automorphism groups of 
the vertex operator algebras for the 
Wess-Zumino-Novikov-Witten models contain the automorphism 
groups of the finite-dimensional Lie algebras that one starts with. 
There are also spectral flow automorphisms of 
the vertex operator algebras for the 
Wess-Zumino-Novikov-Witten models and of superconformal 
vertex operator algebras. Many of these automorphisms
are of infinite orders and do not act semisimply. 
In the case that these automorphisms are the exponentials of 
weight $1$ elements of the vertex operator algebras,
strongly $\C$-graded generalized 
twisted $V$-modules associated to these automorphisms 
are constructed by applying the results of the present paper.

The present paper is organized as follows: In Section 2, 
we reformulate the notion of 
$g$-twisted $V$-module in 
the case of finite-order $g$ using complex variables and 
duality properties. Then we introduce the notions of 
strongly
$\C/\Z$- and $\C$-graded generalized  
$g$-twisted $V$-module in the general case in Section 3. 
In Section 4, we study the exponentials of the
negative-power parts of vertex operators. 
The construction of strongly
$\C$-graded generalized 
$g_{u}$-twisted 
$V$-modules for $u\in V_{(1)}$ is given 
in Section 5. The explicit examples related to 
the triplet $W$-algebras are also given in 
this section.

\paragraph{Acknowledgments}
The author would like to thank Antun Milas for pointing
out that Theorem \ref{main} can be used to show that 
certain intertwining operators constructed in \cite{AM}
are in fact the twisted vertex operator maps for 
certain generalized twisted modules in the sense of
the definition given in the present paper. The author is 
supported in part by NSF grant PHY-0901237.

\renewcommand{\theequation}{\thesection.\arabic{equation}}
\renewcommand{\thethm}{\thesection.\arabic{thm}}
\setcounter{equation}{0}
\setcounter{thm}{0}

\section{Twisted modules for finite-order automorphisms of 
a vertex operator algebra}

The notion
of $g$-twisted module associated to a finite order 
automorphism $g$ of $V$ is formulated by collecting the
properties of twisted vertex operators obtained in \cite{Le1},
\cite{FLM2}, \cite{Le2} and \cite{FLM3}; 
in particular, it uses 
a twisted Jacobi identity (see (\ref{twisted-jacobi}) below) 
as the main axiom. See \cite{FLM3}, \cite{D} and \cite{DL}.
For an infinite-order automorphism $g$ of $V$, 
we cannot write down a generalization of the twisted Jacobi identity
in the case of finite-order automorphisms. The twisted Jacobi identity 
can also be replaced by the weak commutativity and (twisted) weak 
associativity in the case of finite-order automorphisms 
(see (\ref{wk-assoc-2}) below). 
But the side of the weak associativity containing the iterate of 
the twisted vertex operator map
cannot be generalized to the case of an (necessarily infinite-order)
automorphism which acts on the vertex operator algebra or the twisted 
module to be defined nonsemisimply. These difficulties force
us to use the complex variable approach.
To motivate our notion of $g$-twisted module for a general 
automorphism $g$ of $V$, we reformulate the definition of $g$-twisted module 
in terms of complex variables in this section. 

Let $(V, Y, \one, \omega)$ be a vertex operator algebra and 
$g$ an automorphism of $V$ of order $k \in
\mathbb{Z}_+$. Then $V=\coprod_{j\in \Z/k\Z}V^{j}$ where 
$V^{j}=\{v\in V\;|\;gv=\eta^{j}v\}$ for $j\in \Z/k\Z$, 
where $\eta=e^{\frac{2\pi \sqrt{-1}}{k}}$.
We first recall the notion of $g$-twisted $V$-module.

\begin{defn}\label{twisted-mod}
{\rm A {\it $g$-twisted $V$-module} $W$ is a $\C$-graded
vector space $W=\coprod_{n \in \C}W_{(n)}$
equipped with  a linear map
\begin{eqnarray*}
Y^{g}: V \otimes W&\to & W[[x^{1/k},x^{-1/k}]] \nn
v \otimes w&\mapsto & Y^g (v,x)w=\sum_{n\in \frac{1}{k}\Z }
Y_{n}^{g}(v)w x^{-n-1} 
\end{eqnarray*}
satisfying the following conditions:
\begin{enumerate}

\item The {\it grading-restriction condition}:
For each $n\in \C$, $\dim W_{(n)}< \infty $ and $W_{(n+l/k)}=0$ 
for all sufficiently negative integers $l$.

\item The {\it formal monodromy condition}: For 
$j\in \Z/k\Z$ and $v\in V^j$,
$Y^g(v,x)=\sum_{n\in j/k+ \Z}Y_{n}^{g}(v) x^{-n-1}$.

\item The {\it 
lower-truncation condition}: For 
$v\in V$ and $w\in W$,
$Y_{n}^{g}(v) w=0$ for $n$ sufficiently large.

\item The {\it identity property}: $Y^g ({\bf 1},x)=1$.

\item The {\it twisted Jacobi identity}: For $u,v\in V$ and $w\in W$,
\begin{eqnarray}
\lefteqn{x^{-1}_0\delta\left(\frac{x_1-x_2}{x_0}\right)
Y^g(u,x_1)Y^g(v,x_2)}\nn
&&\quad\quad\quad\quad
 -x^{-1}_0\delta\left(\frac{x_2-x_1}{-x_0}\right) Y^g
(v,x_2)Y^g (u,x_1)\nn
&&
= \frac{1}{k}\sum_{j\in \Z /k \Z}x_2^{-1}
\delta\left(\eta^j\frac{(x_1-x_0)^{1/k}}{x_2^{1/k}}\right)Y^g (Y(g^j
u,x_0)v,x_2).\label{twisted-jacobi}
\end{eqnarray}

\item The {\it $L^{g}(0)$-grading condition}: For $n\in \Z$, let
$L^g (n)=Y_{n+1}^g(\omega)$, that, is, 
$Y^g(\omega,x)=\sum_{n\in \Z} L^g (n)x^{-n-2}$. Then 
$L^g (0) w=m w$ for $w \in W_{(m)}$.

\item  The {\it $L^{g}(-1)$-derivative condition}: For $v\in V$,
$$\frac{d}{dx}Y^g (v,x)=Y^g (L(-1)v,x).$$

\end{enumerate}}
\end{defn}

We denote a $g$-twisted $V$-module defined above by
$(W,Y^g)$ (or briefly, by $W$). Note that in the definition 
above, we use $Y_{n}^{g}(u)$ instead of $u_{n}^{g}$
for $n\in \frac{1}{k}\Z$
to denote the components of the vertex operator 
$Y^{g}(u, x)$. In the present paper, we shall always use
such notations to denote the components of a vertex operator.
For example, for our vertex operator algebra $(V, Y, \one, \omega)$,
we shall use $Y_{n}(u)$ instead of $u_{n}$ for $n\in \Z$ to denote 
the components of the vertex operator $Y(u, x)$. 
Similarly, for a $V$-module $(W, Y_{W})$, 
we shall use $(Y_{W})_{n}(u)$ instead of $u_{n}$ for $n\in \Z$ to denote 
the components of the vertex operator $Y_{W}(u, x)$.

This definition of $g$-twisted $V$-module is a natural generalization 
of the definition of $V$-module given in \cite{Bo}, \cite{FLM3} and \cite{FHL}. 
The main axiom is the twisted Jacobi identity. The following 
result is proved by Li \cite{Li}:

\begin{prop}\label{weak-duality}
The twisted Jacobi identity in Definition \ref{twisted-mod}
can be replaced by the following two properties:
\begin{enumerate}

\item The {\it weak commutativity} or {\it locality}:
For $u, v\in V$, there exists $N\in \Z_{+}$ such that
\begin{equation}\label{locality-2}
(x_{1}-x_{2})^{N}Y^{g}(u, x_{1})Y^{g}(v, x_{2})
=(x_{2}-x_{1})^{N}Y^{g}(v, x_{2})Y^{g}(u, x_{1}).
\end{equation}

\item The {\it weak associativity}: For $u\in V^{j}$ and $w\in W$,
there exists $M\in \Z_{+}$ such that 
\begin{equation}\label{wk-assoc-2}
(x_{2}+x_{0})^{\tilde{j}/k+M}Y^{g}(Y(u, x_{0})v, x_{2})
=(x_{0}+x_{2})^{\tilde{j}/k+M}Y^{g}(u, x_{0}+x_{2})Y^{g}(v, x_{2})
\end{equation}
for $v\in V$, where $\tilde{j}\in j$ satisfies 
$0\le \tilde{j}<k$.
\end{enumerate}
\end{prop}

One advantage of the definition, the weak commutativity
and the weak associativity above
is that they use only formal variables so that one 
can discuss the multivaluedness of the correlation functions
algebraically. But as we mentioned above,
this formulation of 
$g$-twisted $V$-module for  $g$ of finite order
cannot be generalized to 
the case that $g$ is of infinite order. So we first reformulate
this definition using complex variables and a duality property.

We need a complex variable formulation of vertex operator maps
and some notations first. 
For any $z\in \C^{\times}$, we shall use $\log z$ to denote 
$\log |z|+\sqrt{-1}\arg z$ where $0\le \arg z<2\pi$. In general, we shall
use $l_{k}(z)$ to denote $\log z+2\sqrt{-1}\pi$ for $k\in \Z$. 
Let $W_{1}$, $W_{2}$ and $W_{3}$ 
be $\C$-graded vector spaces. The $\C$-gradings on $W_{1}$ and $W_{2}$
induce a $\C$-grading on $W_{1}\otimes W_{2}$. 
For a formal variable $x^{1/k}$, we define 
its {\it degree} to be $-1/k$. Then the grading on $W_{3}$ and the degree
of $x$ together give
$W_{3}[[x^{1/k}, x^{-1/k}]]$ a $\C$-graded vector space structure. 
Let
\begin{eqnarray*}
X: W_{1}\otimes W_{2}&\to& W_{3}[[x^{1/k}, x^{-1/k}]]\nn
w_{1}\otimes w_{2}&\mapsto &X(w_{1}, x)w_{2}=\sum_{n\in \frac{1}{k}\Z}
X_{n}(w_{1})w_{2}x^{-n-1}
\end{eqnarray*}
is a linear map preserving the gradings.
For any $w_{1}\in W_{1}$ and $w_{2}\in W_{2}$,
we define $X^{p}(w_{1}, z)w_{2}$ for $p\in \Z /k\Z$
to be the elements 
$X(w_{1}, x)w_{2}\lbar_{x^{n}=e^{nl_{\tilde{p}}(z)}}$
of $\overline{W}_{3}=\prod_{n\in \C}(W_{3})_{(n)}$ 
where for $p\in \Z /k\Z$, 
$\tilde{p}$ is the integer in $p$
satisfying $0\le \tilde{p}<k$. 
When $p=0$, for simplicity, we shall denote $X^{0}(w_{1}, z)w_{2}$
simply as $X(w_{1}, z)w_{2}$.
For $z\in \C^{\times}$, we have the 
linear maps
\begin{eqnarray}
X^{p}(\cdot, z)\cdot: W_{1}\otimes W_{2}&\to& \overline{W}_{3},\nn
w_{1}\otimes w_{2}&\mapsto& X^{p}(w_{1}, z)w_{2}
\end{eqnarray}
and we shall denote $X^{0}(\cdot, z)\cdot$ simply as 
$X(\cdot, z)\cdot$. For any $w_{3}'\in W_{3}'$, 
$$\langle w_{3}', X^{p}(w_{1}, z)w_{2}\rangle$$
are branches of a multivalued analytic function defined on
$\C^{\times}$. We can view this multivalued analytic function 
on $\C^{\times}$ as a single-valued analytic function defined 
on a $k$-fold covering space of $\C^{\times}$. 
For $p\in \Z/k\Z$, we have a map 
\begin{eqnarray*}
X^{p}: \C^{\times}&\to& \hom(W_{1}\otimes W_{2}, \overline{W}_{3})\nn
z&\mapsto& X^{p}(\cdot, z)\cdot.
\end{eqnarray*}
We shall call $X^{p}$ the {\it $p$-th analytic branch}
of the map $X$. 

In particular, for the twisted vertex operator map $Y^{g}$,
we have its $p$-th analytic branch $Y^{g; p}$.
In terms of these analytic branches, we have: 

\begin{prop}\label{equiv-cond}
The formal monodromy condition in Definition 
\ref{twisted-mod} can be replaced
by the following condition: For 
$p\in \Z /k\Z$, $z\in \C^{\times}$ and $v \in V$,
\begin{equation}\label{deck-trans}
Y^{g; p+1} (gv, z) =  Y^{g; p}(v, z).
\end{equation}
\end{prop}
\pf
Let $W$ be a $g$-twisted $V$-module satisfying Definition \ref{twisted-mod}.
Then for $v\in V^{j}$, we have $gv=\eta^{j}v\in V^{j}$.
It is easy to see that this fact and  the formal monodromy condition give
$$Y^{g; p+1}(gv, z)=Y^{g; p}(v, z)$$
for $p\in \Z /k\Z$, 
$z\in \C^{\times}$ and $v\in V^{j}$, 

Conversely, assume that $W$ satisfies the modification of 
Definition \ref{twisted-mod} by replacing the formal monodromy condition
with the equivariance condition. Then it is easy to see that the 
fact $gv=\eta^{j}v\in V^{j}$
for $v\in V^{j}$ and the equivariance condition give
$$\sum_{n\in \frac{1}{k}\Z}Y_{n}^{g}(v)\eta^{j}e^{(-n-1)l_{p+1}(z)}
=\sum_{n\in \frac{1}{k}\Z}Y_{n}^{g}(v)e^{\frac{2\pi (n+1)\sqrt{-1}}{k}}
e^{(-n-1)l_{p+1}(z)}$$
for $j\in \Z/k\Z$ and $v\in V^{j}$.
So we have
$$\sum_{n\in \frac{1}{k}\Z}Y_{n}^{g}(v)(\eta^{j}-
e^{\frac{2\pi (-n-1)\sqrt{-1}}{k}})e^{(n+1)l_{p}(z)}=0,$$
which implies that 
$$Y_{n}^{g}(v)(\eta^{j}-
e^{\frac{2\pi (n+1)\sqrt{-1}}{k}})=0$$
for $n\in \frac{1}{k}\Z$. 
So we have either $Y_{n}^{g}(v)=0$ or 
$\eta^{j}=
e^{\frac{2\pi (n+1)\sqrt{-1}}{k}}$, that is, either $Y_{n}^{g}(v)= 0$
or $n\in j/k+\Z$. So we obtain
$Y^{g}(v, x)=\sum_{n\in j/k+\Z}Y^{g}_{n}(v)x^{-n-1}$.
\epfv

We shall call (\ref{deck-trans}) in Proposition \ref{equiv-cond}
the {\it equivariance property}. 

For $n\in \C$, let $\pi_{n}: W\to W_{(n)}$ be the projection.
We have:

\begin{thm}\label{duality-prop} 
The twisted Jacobi identity in Definition 
\ref{twisted-mod} can be replaced
by the following property:
For any $u, v\in V$,
$w\in W$ and $w'\in W'$, there exists a multivalued analytic function
of the form 
$$f(z_{1}, z_{2})=\sum_{r, s=N_{1}}^{N_{2}}
a_{rs}z_{1}^{r/k}z_{2}^{s/k}(z_{1}-z_{2})^{-N}$$
for $N_{1}, N_{2}\in \Z$ and $N\in \Z_{+}$, such that
the series 
\begin{eqnarray}
&{\displaystyle \langle w', Y^{g; p}(u, z_{1})Y^{g;p}(v, z_{2})w\rangle
=\sum_{n\in \C}\langle w', Y^{g;p}(u, z_{1})\pi_{n}Y^{g;p}(v, z_{2})w\rangle,}&
\label{prod-1-branch}\\
&{\displaystyle \langle w', Y^{g;p}(v, z_{2})Y^{g;p}(u, z_{1})w\rangle
=\sum_{n\in \C}\langle w', Y^{g;p}(v, z_{2})\pi_{n}Y^{g;p}(u, z_{1})w\rangle,}&
\label{prod-2-branch}\\
&{\displaystyle \langle w', Y^{g; p}(Y(u, z_{1}-z_{2})v, z_{2})w\rangle
=\sum_{n\in \C}\langle w', Y^{g; p}(\pi_{n}Y(u, z_{1}-z_{2})v, z_{2})w\rangle}&
\label{iterate-1-branch}\nn
\end{eqnarray}
are absolutely convergent in the regions $|z_{1}|>|z_{2}|>0$,
$|z_{2}|>|z_{1}|>0$, $|z_{2}|>|z_{1}-z_{2}|>0$, respectively,
to the branch 
\begin{equation}\label{correl-fn-branch}
\sum_{r, s=N_{1}}^{N_{2}}
a_{rs}e^{(r/k)l_{p}(z_{1})}e^{(s/k)l_{p}(z_{2})}(z_{1}-z_{2})^{-N}
\end{equation}
of $f(z_{1}, z_{2})$.
\end{thm}
\pf
Assume that $W$ is a $g$-twisted $V$-module satisfying 
Definition 
\ref{twisted-mod}. By Proposition \ref{weak-duality},
there exists $N\in \Z_{+}$ such that
(\ref{locality-2}) holds.

For $w\in W$ and $w'\in W'$, we know that 
\begin{equation}\label{prod-1}
\langle w', Y^{g}(u, x_{1})Y^{g}(v, x_{2})w\rangle
\end{equation}
involves finitely many positive powers of $x_{1}^{1/k}$ and 
finitely many negative powers of $x_{2}^{1/k}$ and
\begin{equation}\label{prod-2}
\langle w', Y^{g}(v, x_{2})Y^{g}(u, x_{1})w\rangle
\end{equation}
involves finitely many negative powers of $x_{1}^{1/k}$ and 
finitely many positive powers of $x_{2}^{1/k}$. Since $(x_{1}-x_{2})^{N}$
is a polynomial in $x_{1}^{1/k}$ and $x_{2}^{1/k}$, the same
statements are also true for 
$$(x_{1}-x_{2})^{N}\langle w', Y^{g}(u, x_{1})Y^{g}(v, x_{2})w\rangle$$
and 
$$(x_{1}-x_{2})^{N}\langle w', Y^{g}(v, x_{2})Y^{g}(u, x_{1})w\rangle,$$
respectively. By (\ref{locality-2}), they are actually equal. 
Thus they are both equal to
$$\sum_{r, s=N_{1}}^{N_{2}}
a_{rs}x_{1}^{r/k}x_{2}^{s/k}\in \C[x_{1}^{1/k}, x_{1}^{-1/k}, 
x_{2}^{1/k}, x_{2}^{-1/k}].$$
Since (\ref{prod-1})
involves finitely many positive powers of $x_{1}^{1/k}$ and 
finitely many negative powers of $x_{2}^{1/k}$,
it must be equal to the expansion of 
$$\sum_{r, s=N_{1}}^{N_{2}}
a_{rs}x_{1}^{r/k}x_{2}^{s/k}(x_{1}-x_{2})^{-N}$$
as a series of the same form. Thus the expansion of (\ref{correl-fn-branch})
in the region $|z_{1}|>|z_{2}|>0$ is equal to 
(\ref{prod-1-branch}). The same argument shows that 
the expansion of (\ref{correl-fn-branch})
in the region $|z_{2}|>|z_{1}|>0$ is equal to 
(\ref{prod-2-branch}).

For $j\in \Z/k\Z$, let $\tilde{j}\in j$ satisfying $0\le \tilde{j}<k$.
For $u\in V^{j}$, by Proposition \ref{weak-duality}
there exists $M\in \Z_{+}$ such that (\ref{wk-assoc-2})
holds.

For $w\in W$ and $w'\in W'$, we know that (\ref{prod-1})
involves finitely many positive powers of $x_{1}^{1/k}$ and 
finitely many negative powers of $x_{2}^{1/k}$. So
\begin{equation}\label{prod-1-1}
\langle w', Y^{g}(u, x_{0}+x_{2})Y^{g}(v, x_{2})w\rangle
\end{equation}
involves finitely many positive powers of $x_{0}^{1/k}$ and 
finitely many negative powers of $x_{2}^{1/k}$. We also know that
\begin{equation}\label{iter}
\langle w', Y^g(Y(u, x_0)v, x_2)w\rangle
\end{equation}
involves finitely many negative powers of $x_{0}^{1/k}$ and 
finitely many positive powers of $x_{2}^{1/k}$. Since 
$(x_{0}+x_{2})^{\tilde{j}/k+M}$
involves finitely many positive powers of $x_{0}^{1/k}$ and 
no negative powers of $x_{2}^{1/k}$ and 
$(x_{2}+x_{0})^{\tilde{j}/k+M}$ 
involves finitely many positive powers of $x_{2}^{1/k}$ and 
no negative powers of $x_{2}^{1/k}$,  the same
statements above for (\ref{prod-1-1}) and (\ref{iter}) 
are also true for 
$$(x_{0}+x_{2})^{\tilde{j}/k+M}
\langle w', Y^{g}(u, x_{0}+x_{2})Y^{g}(v, x_{2})w\rangle$$
and 
$$(x_{2}+x_{0})^{\tilde{j}/k+M}\langle w', Y^g(Y(u, x_0)v, x_2)w\rangle,$$
respectively. 
Since we have proved that (\ref{prod-1}) is equal to 
$$\sum_{r, s=N_{1}}^{N_{2}}
a_{rs}x_{1}^{r/k}x_{2}^{s/k}(x_{1}-x_{2})^{-N},$$
we obtain
\begin{eqnarray*}
\lefteqn{(x_{0}+x_{2})^{\tilde{j}/k+M}
\langle w', Y^{g}(u, x_{0}+x_{2})Y^{g}(v, x_{2})w\rangle}\nn
&&=\sum_{r, s=N_{1}}^{N_{2}}
a_{rs}(x_{0}+x_{2})^{r/k+\tilde{j}/k+M}x_{2}^{s/k}x_{0}^{-N}.
\end{eqnarray*}
By (\ref{wk-assoc-2}), we also have
\begin{eqnarray*}
\lefteqn{(x_{2}+x_{0})^{\tilde{j}/k+M}\langle w', Y^g(Y(u, x_0)v, x_2)w\rangle}\nn
&&=\sum_{r, s=N_{1}}^{N_{2}}
a_{rs}(x_{2}+x_{0})^{r/k+\tilde{j}/k+M}x_{2}^{s/k}x_{0}^{-N}.
\end{eqnarray*}
Thus the expansion of (\ref{correl-fn-branch})
in the region $|z_{2}|>|z_{1}-z_{2}|>0$ is equal to 
(\ref{iterate-1-branch}). 

Conversely, assume that the property in the theorem holds. Let 
$$f_{j}(x_{1}, x_{2}, x_{0})=\sum_{r, s=N_{1}}^{N_{2}}a_{rs}\eta^{-jr}
x_{1}^{r/k}x_{2}^{s/k}x_{0}^{-N}$$
for $j\in \Z/kZ$.
Then this property gives
\begin{eqnarray}
\langle w', Y^{g}(u, x_{1})Y^{g}(v, x_{2})w\rangle
&=&f_{0}(x_{1}, x_{2}, x_{1}-x_{2}),\label{formal-prod-1}\\
\langle w', Y^{g}(v, x_{2})Y^{g}(u, x_{1})w\rangle
&=&f_{0}(x_{1}, x_{2}, -x_{2}+x_{1}), \label{formal-prod-2}\\
\langle w', Y^{g}(Y(u, z_{1}-z_{2})v, z_{2})w\rangle
&=&f_{0}(x_{2}+x_{0}, x_{2}, x_{0}) \label{formal-iterate-1}
\end{eqnarray}
for $u, v\in V$, $w\in W$ and $w'\in W'$.
{}From (\ref{formal-prod-1}), 
(\ref{formal-iterate-1}) and the equivariance property,
we have
\begin{equation}\label{formal-iterate-2}
\langle w', Y^{g}(Y(g^{j}u, z_{1}-z_{2})v, z_{2})w\rangle
=f_{j}(x_{2}+x_{0}, x_{2}, x_{0})
\end{equation}
for $u, v\in V$, $w\in W$ and $w'\in W'$.

We have the formal variable identity
$$x_{0}^{-1}\delta\left(\frac{x_{1}-x_{2}}{x_{0}}\right)
-x_{0}^{-1}\delta\left(\frac{x_{2}-x_{1}}{-x_{0}}\right)
=\frac{1}{k}\sum_{j\in \Z/k\Z}x_{2}^{-1}
\delta\left(\eta^{j}\frac{(x_{1}-x_{0})^{1/k}}{x_{2}^{1/k}}\right).$$
Using this identity, the properties of the delta-function and 
the explicit form of $f_{0}(x_{1}, x_{2}, x_{0})$, we can prove
\begin{eqnarray}\label{formal-calc}
&{\displaystyle x_{0}^{-1}\delta\left(\frac{x_{1}-x_{2}}{x_{0}}\right)
f_{0}(x_{1}, x_{2}, x_{1}-x_{2})
-x_{0}^{-1}\delta\left(\frac{x_{2}-x_{1}}{-x_{0}}\right)
f_{0}(x_{1}, x_{2}, -x_{2}+x_{1})}&\nn
&{\displaystyle =\frac{1}{k}\sum_{j\in \Z/k\Z}x_{2}^{-1}
\delta\left(\eta^{j}\frac{(x_{1}-x_{0})^{1/k}}{x_{2}^{1/k}}\right)
f_{j}(x_{2}+x_{0}, x_{2}, x_{0}).}&
\end{eqnarray}
{}From (\ref{formal-prod-1})--(\ref{formal-iterate-1}) 
and (\ref{formal-calc}),
we obtain 
\begin{eqnarray*}
\lefteqn{x^{-1}_0\delta\left(\frac{x_1-x_2}{x_0}\right)
\langle w', Y^g(u,x_1)Y^g(v,x_2)w\rangle}\nn
&&\quad\quad\quad\quad
 -x^{-1}_0\delta\left(\frac{x_2-x_1}{-x_0}\right) 
\langle w', Y^g(v,x_2)Y^g (u,x_1)w\rangle\nn
&&=\frac{1}{k}\sum_{j\in \Z /k \Z} x_2^{-1}
\delta\left(\eta^j\frac{(x_1-x_0)^{1/k}}{x_2^{1/k}}\right)
\langle w', Y^g (Y(g^j
u,x_0)v,x_2)w\rangle
\end{eqnarray*}
for $u, v\in V$, $w\in W$ and $w'\in W'$.
Since $w\in W$ and $w'\in W'$ are arbitrary,
we obtain the twisted Jacobi identity (\ref{twisted-jacobi})
for $u, v\in V$.
\epfv

We shall call the property in Theorem \ref{duality-prop} 
the {\it duality property}.

\begin{rema}
{\rm For simplicity, 
we give only the definition of twisted 
module for a vertex operator algebra, not for a vertex operator 
superalgebra. But it is trivial to generalize it
to that of twisted module for a vertex 
operator superalgebra. }
\end{rema}

\renewcommand{\theequation}{\thesection.\arabic{equation}}
\renewcommand{\thethm}{\thesection.\arabic{thm}}
\setcounter{equation}{0}
\setcounter{thm}{0}

\section{Definitions of strongly $\C/\Z$- and $\C$-graded generalized 
$g$-twisted module}

We give a definition of strongly $\C/\Z$-graded generalized 
$g$-twisted $V$-module
for an automorphism $g$ of $V$ of possibly infinite order in this section. 
In the case that $V$ admits an additional $\C$-grading compatible with
$g$, we also give a definition of strongly $\C$-graded generalized 
$g$-twisted $V$-module.
In the case that $g$ is of finite order, $g$ always acts on 
$V$ semisimply so that we have $V=\coprod_{j\in \Z/k\Z}V^{j}$.
When the order of $g$ is infinite, it does not have to act on 
$V$ semisimply. In particular, we have to allow not only nonintegral 
powers of $x$ 
in the twisted vertex operators, but also integral powers of $\log x$
in these operators. As we have mentioned in the introduction, 
the left-hand side of the twisted Jacobi identity cannot be 
generalized to the case that the order of $g$ is infinite. 
Also we cannot straightforwardly 
generalize the formal monodromy condition to the general case. 
On the other hand, the equivariance property and the duality 
property in Proposition \ref{equiv-cond} and Theorem \ref{duality-prop},
respectively, can be 
generalized easily. Our definition will use the generalizations of these
axioms.

First we need to further 
generalize our consideration of analytic branches of 
the map $X$ in the preceding section. For formal variables 
$x$ and $\log x$, we define 
their {\it degrees} to be $-1$ and $0$. Let $W_{1}$, $W_{2}$ and $W_{3}$ 
be $\C$-graded vector spaces. Then the 
grading on $W_{3}$ and the degrees
of $x$ and $\log x$ together give
$W_{3}\{x\}[\log x]$ a $\C$-grading.
Let
\begin{eqnarray*}
X: W_{1}\otimes W_{2}&\to& W_{3}\{x\}[\log x]\nn
w_{1}\otimes w_{2}&\mapsto &X(w_{1}, x)w_{2}=\sum_{n\in \C}
\sum_{k=1}^{K}
X_{n, k}(w_{1})w_{2}x^{-n-1}(\log x)^{k}
\end{eqnarray*}
is a linear map preserving the gradings.
For any $w_{1}\in W_{1}$ and $w_{2}\in W_{2}$,
we define $X^{p}(w_{1}, z)w_{2}$ for $p\in \Z$
to be the elements 
$X(w_{1}, x)w_{2}\lbar_{x^{n}=e^{nl_{p}(z)},\;\log x=l_{p}(z)}$
of $\overline{W}_{3}$. 
When $p=0$, we shall denote $X^{0}(w_{1}, z)w_{2}$
simply as $X(w_{1}, z)w_{2}$.
For $z\in \C^{\times}$, we have the 
linear maps
\begin{eqnarray}
X^{p}(\cdot, z)\cdot: W_{1}\otimes W_{2}&\to& \overline{W}_{3},\nn
w_{1}\otimes w_{2}&\mapsto& X^{p}(w_{1}, z)w_{2}
\end{eqnarray}
and we shall denote $X^{0}(\cdot, z)\cdot$ simply as 
$X(\cdot, z)\cdot$. For any $w_{3}'\in W_{3}'$, 
$$\langle w_{3}', X^{p}(w_{1}, z)w_{2}\rangle$$
are branches of a multivalued analytic function defined on
$\C^{\times}$. We can view this multivalued analytic function 
on $\C^{\times}$ as a single-valued analytic function defined 
on a  covering space of $\C^{\times}$. 
For $p\in \Z$, we have a map 
\begin{eqnarray*}
X^{p}: \C^{\times}&\to& \hom(W_{1}\otimes W_{2}, \overline{W}_{3})\nn
z&\mapsto& X^{p}(\cdot, z)\cdot.
\end{eqnarray*}
We shall call $X^{p}$ the {\it $p$-th analytic branch}
of the map $X$. 

We also need a $\C^{\times}$-, or equivalently, $\C/\Z$-grading 
structure on $V$ given by an automorphism 
$g$ of $V$. Since $g$ preserve $V_{(n)}$ for $n\in \Z$
and $V_{(n)}$ for $n\in \Z$ are all finite dimensional, we have
$$V_{(n)}=\coprod_{a\in \C^{\times}}V_{(n)}^{[a]}$$
for $n\in \Z$, where for $n\in \Z$ and $a\in \C^{\times}$,
if $a$ is an eigenvalue for the operator $g$ on $V_{(n)}$,
$V_{(n)}^{[a]}$ is the generalized eigenspace of $V_{(n)}$ 
for $g$ with the eigenvalue $a$, and if $a$ is not an eigenvalue 
for $g$ on $V_{(n)}$, $V_{(n)}^{[a]}=0$. Since the multiplicative
abelian group $\C^{\times}$ is isomorphic to the additive abelian group 
$\C/\Z$ through the isomorphism $a\mapsto 
\frac{1}{2\pi \sqrt{-1}}\log a + \Z$, we can use $\C/\Z$
instead of $\C^{\times}$ so that 
$$V_{(n)}=\coprod_{\alpha\in \C/\Z}V_{(n)}^{[\alpha]}$$
and $V$ is doubly graded
by $\C$ and $\C/\Z$ as 
$$V=\coprod_{n\in \Z, \alpha\in \C/\Z}V_{(n)}^{[\alpha]}.$$
Let $V^{[\alpha]}=\coprod_{n\in \Z}V_{(n)}^{[\alpha]}$ for 
$\alpha\in \C/\Z$.
Then 
$V=\coprod_{\alpha\in \C/\Z}V^{[\alpha]}$
and for $v\in V^{[\alpha]}$, $gv=e^{2\pi \sqrt{-1}\;\alpha}v$.
It is easy to see that the vertex operator algebra
$V$ together with this $\C/\Z$-grading is a strongly 
$\C/\Z$-graded vertex operator algebra in the sense
of \cite{HLZ2}. 

If, in addition to the $\C$-grading given by the eigenvalues of $L(0)$,
$V$ has another $\C$-grading $V=\coprod_{\alpha\in \C}V^{[\alpha]}$
(not just a $\C/\Z$-grading as we have given above)
which is compatible with $g$ (that is,
for $v\in V^{[\alpha]}$,  $gv=e^{2\pi \sqrt{-1}\;\alpha}v$), 
then $V=\coprod_{n, \alpha\in \C}V_{(n)}^{[\alpha]}$
equipped with this $\C$-grading has a strongly 
$\C$-graded vertex operator algebra structure.

For a strongly 
$\C/\Z$-graded (or $\C$-graded) vertex operator algebra, a natural 
module category is the category of strongly $\C/\Z$-graded (or $\C$-graded,
respectively)
generalized $V$-modules. Moreover, if the $\C/\Z$-grading (or $\C$-graded)
is given by (or compatible with) $g$, it is natural to 
consider  strongly  $\C/\Z$-graded (or $\C$-graded) generalized $V$-modules
with actions of $g$ such that the $\C^{\times}$-gradings 
(or $\C$-gradings) are given 
by (or compatible with) generalized eigenspaces for $g$ in the obvious sense.
The category of the generalizations of $g$-twisted modules should be 
a generalization of the category of such strongly $\C/\Z$- or $\C$-graded
generalized $V$-modules. In particular, there should also
be an action of $g$
and also a $\C/\Z$-grading (or $\C$-grading) given by 
(or compatible with) generalized eigenspaces for 
the action of $g$. 
Here is our generalization of 
$g$-twisted $V$-module for finite-order $g$ to the general case:

\begin{defn}
{\rm Let $(V,Y,{\bf 1},\omega)$ be a vertex operator algebra
and let $g$ be an automorphism of $V$. 
A {\it strongly $\C/\Z$-graded generalized 
$g$-twisted $V$-module} is a $\C\times \C/\Z$-graded
vector space
$W=\coprod_{n\in \C, \alpha\in \C/\Z}W_{[n]}^{[\alpha]}$
equipped with a linear map
\begin{eqnarray*}
Y^{g}: V \otimes W&\to & W\{x\}[\log x] \\
v \otimes w&\mapsto & Y^g (v, x)w
\end{eqnarray*}
and an action of $g$
satisfying the following conditions:

\begin{enumerate}

\item The {\it grading-restriction condition}: 
For each $n\in \C, \alpha\in \C/\Z$, 
$\dim W_{[n]}^{[\alpha]}< \infty$ and 
$W_{[n+l]}^{[\alpha]}=0$ 
for all sufficiently negative real number $l$.

\item The {\it equivariance property}:
For $p\in \Z$, $z\in \C^{\times}$, $v\in V$ and $w\in W$,
$Y^{g; p+1} (gv, z)w =  Y^{g; p}(v, z)w$,
where for $p\in \Z$, $Y^{g;p}$ is the $p$-th analytic branch of 
$Y^{g}$.

\item The {\it identity property}: For $w\in W$, 
$Y^g ({\bf 1},x)w=w.$

\item The {\it duality property}: Let $W'=\coprod_{n\in \C, 
\alpha\in \C/\Z}
(W_{[n]}^{[\alpha]})^{*}$ and, for $n\in \C$,
$\pi_{n}: W\to W_{[n]}$ be the projection. For any $u, v\in V$,
$w\in W$ and $w'\in W'$, there exists a multivalued analytic function
of the form 
$$f(z_{1}, z_{2})=\sum_{i, j, k, l=0}^{N}a_{ijkl}z_{1}^{m_{i}}
z_{2}^{n_{i}}(\log z_{1})^{k}(\log z_{2})^{l}(z_{1}-z_{2})^{-t}$$
for $N\in \N$, $m_{1}, \dots, m_{N}, n_{1}, \dots, n_{N}\in 
\C$ and $t\in \Z_{+}$, such that
the series 
\begin{eqnarray}
&{\displaystyle \langle w', Y^{g; p}(u, z_{1})Y^{g;p}(v, z_{2})w\rangle
=\sum_{n\in \C}\langle w', Y^{g;p}(u, z_{1})\pi_{n}Y^{g;p}(v, z_{2})w\rangle,}
\quad&\\
&{\displaystyle \langle w', Y^{g;p}(v, z_{2})Y^{g;p}(u, z_{1})w\rangle
=\sum_{n\in \C}\langle w', Y^{g;p}(v, z_{2})\pi_{n}Y^{g;p}(u, z_{1})w\rangle,}
\quad&\\
&{\displaystyle \langle w', Y^{g; p}(Y(u, z_{1}-z_{2})v, z_{2})w\rangle
=\sum_{n\in \C}\langle w', Y^{g; p}(\pi_{n}Y(u, z_{1}-z_{2})v, z_{2})w\rangle}&\nn
\end{eqnarray}
are absolutely convergent in the regions $|z_{1}|>|z_{2}|>0$,
$|z_{2}|>|z_{1}|>0$, $|z_{2}|>|z_{1}-z_{2}|>0$, respectively,
to the branch 
$$\sum_{i, j, k, l=0}^{N}
a_{ijkl}e^{m_{i}l_{p}(z_{1})}e^{n_{j}l_{p(z_{2})}}
l_{p}(z_{1})^{k}l_{p}(z_{2})^{l}(z_{1}-z_{2})^{-t}$$
of $f(z_{1}, z_{2})$.

\item The {\it $L^{g}(0)$- and $g$-grading conditions}: Let $L^{g}(n)$
for $n\in \Z$ be the operator on $W$ given by $Y^{g}(\omega, x)
=\sum_{n\in \Z}L^{g}(n)^{-n-2}$. Then for $n\in \C$ and $\alpha\in \C/\Z$,
$w \in W_{[n]}^{[\alpha]}$,
there exists $K, \Lambda\in \Z_{+}$ such that $(L^g (0)-n)^{K} w
=(g-e^{2\pi \sqrt{-1}\;\alpha})^{\Lambda}w=0$.

\item The {\it $L(-1)$-derivative condition}: For $v\in V$,
$$\frac{d}{dx}Y^g (v,x)=Y^g (L(-1)v,x).$$
 \end{enumerate}
If $V$ has a strongly $\C$-graded vertex operator algebra
structure compatible with $g$, then a {\it strongly $\C$-graded
generalized $g$-twisted $V$-module} is a $\C\times \C$-graded
vector space
$W=\coprod_{n, \alpha\in \C}W_{[n]}^{[\alpha]}$
equipped with a linear map
\begin{eqnarray*}
Y^{g}: V \otimes W&\to & W\{x\}[\log x] \\
v \otimes w&\mapsto & Y^g (v, x)w
\end{eqnarray*}
and an action of $g$, satisfying the same axioms 
above except that $\C/\Z$ is replaced by $\C$ and the 
$g$-grading condition is replaced by the following 
{\it grading-compatibility condition}:
For $\alpha, \beta\in \C$, $v\in V^{[\alpha]}$ and $w\in W^{[\beta]}$,
$Y^{g}(v, x)w\in W^{[\alpha+\beta]}\{x\}[\log x]$.}
\end{defn}

\begin{rema}
{\rm As in the preceding section, 
for simplicity,  we give only the definition of such generalized
twisted 
module for a vertex operator algebra, not for a vertex operator 
superalgebra. It is also trivial in this general case
to generalize the definition above
to that of twisted module for a vertex 
operator superalgebra. }
\end{rema}

\renewcommand{\theequation}{\thesection.\arabic{equation}}
\renewcommand{\thethm}{\thesection.\arabic{thm}}
\setcounter{equation}{0}
\setcounter{thm}{0}

\section{Exponentiating integrals of negative parts of vertex operators}

In this section we shall exponentiate integrals of negative parts
of vertex operators and prove formulas needed in the next section.
We need first several commutator formulas.

The first is the commutator formula for vertex operators which 
is the special case of the Jacobi identity when we take 
the coefficients of $x_{0}^{-1}$ and expand the formal delta-function
on the right-hand side explicitly:

\begin{prop}
Let $(V, Y, \mathbf{1}, \omega)$ be a
vertex operator algebra and $(W, Y_{W})$ a weak $V$-module. 
The for $u, v\in V$, we have the commutator formula 
\begin{equation}\label{Delta-1}
[Y_{W}(u, x_{1}), Y_{W}(v, x_{2})]
=\sum_{m\in \Z}\sum_{k\in \N}{m\choose k}x_{1}^{-m-1}x_{2}^{m-k}
Y_{W}(Y_{k}(u)v, x_{2}).
\end{equation}
\end{prop}

Let $u\in V$. In this section,  we do not require that the conformal 
weight of $u$ be $1$. 
Let $Y_{W}^{\le -2}(u, x)$ be the part of $Y_{W}(u, x)$ consisting
all the monomials in $x$ with powers of $x$ less than or equal to 
$-2$, that is, 
$$Y_{W}^{\le -2}(u, x)=\sum_{m\in \Z_{+}}Y_{m}(u)x^{-m-1}.$$

For any vector space $W$ and 
$$f(x)=\sum_{n\in \Z\setminus \{-1\}}a_{n}x^{n}
\in W[[x]]+x^{-2}W[[x^{-1}]],$$
let 
$$\int_{0}^{x} f(y)=\sum_{n\in \Z\setminus \{-1\}}\frac{a_{n}}{n+1}x^{n+1}.$$
Then we obtain a linear map ({\it integration from $0$ to $x$}) 
$$\int_{0}^{x}: W[[x]]+x^{-2}W[[x^{-1}]]
\to xW[[x]]+x^{-1}W[[x^{-1}]].$$
We have:

\begin{prop}
For $u, v\in V$,
\begin{eqnarray}\label{Delta-3}
\lefteqn{\left[ \left(\int_{0}^{x_{1}} Y_{W}^{\le -2}(u, y)\right), 
Y_{W}(v, x_{2})\right]}\nn
&&=
Y_{W}\left(\int_{0}^{-x_{1}-x_{2}} Y^{\le -2}(u, y)v, 
x_{2}\right)
+\log\left(1+\frac{x_{2}}{x_{1}}\right)
Y_{W}(Y_{0}(u)v, x_{2}),\nn
\end{eqnarray}
where $\int_{0}^{-x_{1}-x_{2}} Y^{\le -2}(u, y)v$ is defined 
by expanding the series in powers of $-x_{1}-x_{2}$ 
as a series in positive powers of $x_{2}$ and in powers of $x_{1}$.
\end{prop}
\pf
{}From (\ref{Delta-1}), we obtain 
\begin{eqnarray}\label{Delta-2}
[Y_{W}^{\le -2}(u, x_{1}), Y_{W}(v, x_{2})]
&=&\sum_{m\in \Z_{+}}\sum_{k\in \N}{m\choose k}x_{1}^{-m-1}x_{2}^{m-k}
Y_{W}(Y_{k}(u)v, x_{2})\nn
&=&\sum_{k, l\in \N, k+l\ne 0}{k+l\choose k}x_{1}^{-k-l-1}x_{2}^{l}
Y_{W}(Y_{k}(u)v, x_{2}).\nn
\end{eqnarray}
Then 
\begin{eqnarray*}
\lefteqn{\left[\int_{0}^{-x_{1}} Y_{W}^{\le -2}(u, y), 
Y_{W}(v, x_{2})\right]}\nn
&&=\sum_{k, l\in \N, k+l\ne 0}{k+l\choose k}\frac{(-1)^{k+l+1}}{k+l}
x_{1}^{-k-l}x_{2}^{l}
Y_{W}(Y_{k}(u)v, x_{2})\nn
&&=\sum_{k\in \Z_{+}}\frac{(-1)^{k}}{-k}\sum_{l\in \N}
{-k\choose l}x_{1}^{-k-l}x_{2}^{l}
Y_{W}(Y_{k}(u)v, x_{2})\nn
&&\quad +\sum_{l\in \Z_{+}}\frac{(-1)^{l+1}}{l}x_{1}^{-l}x_{2}^{l}
Y_{W}(Y_{0}(u)v, x_{2})\nn
&&=\sum_{k\in \Z_{+}}\frac{(-1)^{k}}{-k}(x_{1}+x_{2})^{-k}
Y_{W}(Y_{k}(u)v, x_{2})
+\log\left(1+\frac{x_{2}}{x_{1}}\right)
Y_{W}(Y_{0}(u)v, x_{2})\nn
&&=Y_{W}\left(\int_{0}^{-x_{1}-x_{2}}Y^{\le -2}(u, y)v, x_{2}\right)
+\log\left(1+\frac{x_{2}}{x_{1}}\right)
Y_{W}(Y_{0}(u)v, x_{2}).
\end{eqnarray*}
\epfv

We also have:

\begin{prop}
For $u\in V$,
\begin{equation}\label{Delta-5}
\left[L(-1), \int_{0}^{-x} Y_{W}^{\le -2}(u, y)\right]
=-\frac{d}{dx}\int_{0}^{x} Y_{W}^{\le -2}(u, y)
-Y_{0}(u)x^{-1}.
\end{equation}
\end{prop}
\pf
By the $L(-1)$-derivative property, we have 
$$[L(-1), Y_{W}(u, y)]=\frac{d}{dy}Y_{W}(u, x),$$
which gives
\begin{equation}\label{Delta-4}
[L(-1), Y_{W}^{\le -2}(u, y)]=\frac{d}{dy}Y_{W}^{\le -2}(u, y)
-Y_{0}(u)y^{-2}.
\end{equation}
Note that the integration $\int_{0}^{x}: W[[x]]+x^{-2}W[[x^{-1}]]
\to xW[[x]]+x^{-1}W[[x^{-1}]]$ defined above is invertible and 
its inverse is in fact the restriction of $\frac{d}{dx}$ to 
the space $xW[[x]]+x^{-1}W[[x^{-1}]]$. 
Integrating both sides of (\ref{Delta-4}) with respect to $y$ 
and then substituting $-x$ for $y$, we obtain
\begin{eqnarray*}
\left[L(-1), \int_{0}^{-x} Y_{W}^{\le -2}(u, y)\right]
&=&Y_{W}^{\le -2}(u, -x)-Y_{0}(u)x^{-1}\nn
&=&-\frac{d}{dx}\int_{0}^{-x} Y_{W}^{\le -2}(u, y)-Y_{0}(u)x^{-1}.
\end{eqnarray*}
\epfv

Let
$$\widetilde{\Delta}_{W}^{(u)}(x)
=\exp\left(\left(\int_{0}^{-x} Y_{W}^{\le -2}(u, y)\right)
+Y_{0}(u)\log x\right).$$
Then we see that $\Delta_{W}^{(u)}(x)\in (\mbox{\rm End})[[x^{-1}]]
[[\log x]].$

\begin{prop}
For $u, v\in V$, 
\begin{eqnarray}
\widetilde{\Delta}_{W}^{(u)}(x)Y(v, x_{2})
&=&Y(\widetilde{\Delta}_{W}^{(u)}(x+x_{2})v, x_{2})
\widetilde{\Delta}_{W}^{(u)}(x),
\label{Delta-6}\\
{[L(-1), \widetilde{\Delta}_{W}^{(u)}(x)]}&=&-\frac{d}{dx}
\widetilde{\Delta}_{W}^{(u)}(x).
\label{Delta-7}
\end{eqnarray}
\end{prop}
\pf
By (\ref{Delta-3}) and (\ref{Delta-5}), we obtain 
(\ref{Delta-6}) and (\ref{Delta-7}), respectively.
\epfv

\renewcommand{\theequation}{\thesection.\arabic{equation}}
\renewcommand{\thethm}{\thesection.\arabic{thm}}
\setcounter{equation}{0}
\setcounter{thm}{0}

\section{Generalized twisted $V$-modules associated to weight $1$ elements of $V$}

In this section, we 
construct strongly $\C$-graded 
generalized twisted $V$-modules 
associated
to $u\in V_{(1)}$ under the  assumption that 
$V_{(0)}=\C\mathbf{1}$ and $V_{(n)}=0$
for $n<0$ and
that $L(1)u=0$.

For such an element $u$, we have an operator $Y_{0}(u)=\res_{x} Y(u, x)$
on $V$. Since the conformal weight of $Y_{0}(u)$ is $0$, 
it preserves the finite-dimensional homogeneous subspaces
$V_{(n)}$ for $n\in \C$. In particular, for any $n\in \C$,
$V_{(n)}=\coprod_{\alpha\in \C}V_{(n)}^{[\alpha]}$
where $V_{(n)}^{[\alpha]}$ is the generalized eigenspaces 
for $Y_{0}(u)$ with eigenvalue $\alpha$ if $\alpha$ is 
an eigenvalue of $Y_{0}(u)$ and is $0$ for other $\alpha\in \C$. 
Thus $V$ has an additional $\C$-grading and, equipped with this 
$\C$-grading, is a 
strongly $\C$-graded vertex operator algebra.

Given a strongly $\C$-graded generalized $V$-module, we would like to 
construct a strongly $\C$-graded generalized twisted $V$-module 
by modifying the vertex operator 
map for the module using the operators $\widetilde{\Delta}_{V}^{(u)}(x)$ 
constructed in the preceding section. But in general 
$\widetilde{\Delta}_{V}^{(u)}(x)v\in V[[x^{-1}]][[\log x]]$ for $v\in V$.
On the other hand, in our definition 
of strongly $\C^{\times}$-graded generalized 
twisted module, the image of the vertex operator map must be in 
$V\{x\}[\log x]$.
So we cannot use $\widetilde{\Delta}_{V}^{(u)}(x)$ directly.
Instead, we first need to construct a series $\Delta_{V}^{(u)}(x)$ 
by modifying $\widetilde{\Delta}_{V}^{(u)}(x)$ and
discuss under what assumptions, 
$\Delta_{V}^{(u)}(x)v$ 
is in $V\{x\}[\log x]$ for $v\in V$.

Let $(Y_{0}(u))_{ss}$ be the semisimple part of $Y_{0}(u)$, that is,
for any generalized eigenvector $v\in V$ for $Y_{0}(u)$ with eigenvalue
$\lambda$, $(Y_{0}(u))_{ss}v=\lambda v$. Then 
$Y_{0}(u)=(Y_{0}(u))_{ss}+(Y_{0}(u)-(Y_{0}(u))_{ss})$ and 
$(Y_{0}(u)-(Y_{0}(u))_{ss})^{N_{u}}=0$. 
The commutator formula (\ref{Delta-1}) for vertex operators gives
\begin{equation}\label{comm-u0}
[Y_{0}(u), Y(v, x_{2})]
=Y(Y_{0}(u)v, x_{2})
\end{equation}
for any $v\in V$. 

The commutator formula (\ref{comm-u0}) gives 
\begin{eqnarray}
{[(Y_{0}(u))_{ss}, Y(v, x_{2})]}
&=&Y((Y_{0}(u))_{ss}v, x_{2}),\label{comm-u0-ss}\\
{[(Y_{0}(u)-(Y_{0}(u))_{ss}), Y(v, x_{2})]}
&=&Y((Y_{0}(u)-(Y_{0}(u))_{ss})v, x_{2}).\label{comm-u0-nil}
\end{eqnarray}
The statement that (\ref{comm-u0-ss}) is a consequence of (\ref{comm-u0})
is a special case of 
a very general fact in linear algebra. Here we give 
a direct proof in this special case. 
Let $v, w\in V$ be generalized eigenvectors for $Y_{0}(u)$ with 
eigenvalues $\lambda_{1}$ and $\lambda_{2}$, respectively. Then there
exist $N_{1}, N_{2}\in \Z_{+}$ such that 
$(Y_{0}(u)-\lambda_{1})^{N_{1}}v=(Y_{0}(u)-\lambda_{2})^{N_{1}}w=0$.
From (\ref{comm-u0}), we obtain
\begin{eqnarray*}
\lefteqn{(Y_{0}(u)-\lambda_{1}-\lambda_{2})Y(v, x_{2})w}\nn
&&=Y((Y_{0}(u)-\lambda_{1})v, x_{2})w+
Y(v, x_{2})(Y_{0}(u)-\lambda_{2})w.
\end{eqnarray*}
Using this formula repeatedly, we obtain
\begin{eqnarray*}
\lefteqn{(Y_{0}(u)-\lambda_{1}-\lambda_{2})^{N_{1}+N_{2}}Y(v, x_{2})w}\nn
&&=\sum_{i=0}^{N_{1}+N_{2}}{N_{1}+N_{2}\choose i}
Y((Y_{0}(u)-\lambda_{1})^{i}v, x_{2})(Y_{0}(u)-\lambda_{2})^{N_{1}+N_{2}-i}w\nn
&&=0.
\end{eqnarray*}
This shows that $Y(v, x_{2})w$ is a formal series whose coefficients 
are generalized eigenvectors for $Y_{0}(u)$ with 
eigenvalue $\lambda_{1}+\lambda_{2}$. Thus
\begin{eqnarray*}
(Y_{0}(u))_{ss}Y(v, x_{2})w&=&(\lambda_{1}+\lambda_{2})Y(v, x_{2})w\nn
&=&Y((Y_{0}(u))_{ss}v, x_{2})w+Y(v, x_{2})(Y_{0}(u))_{ss}w,
\end{eqnarray*}
which is equivalent to 
$$[(Y_{0}(u))_{ss}, Y(v, x_{2})]w
=Y((Y_{0}(u))_{ss}v, x_{2})w.$$
Since $v$ and $w$ are arbitrary generalized eigenvectors for $Y_{0}(u)$,
we obtain (\ref{comm-u0-ss}). The formula (\ref{comm-u0-nil}) 
follows immediately from (\ref{comm-u0}) and (\ref{comm-u0-ss}).

For any
generalized eigenvector $v\in V$ for $Y_{0}(u)$ with eigenvalue
$\lambda$, we have $e^{tY_{0}(u)}v=e^{t(Y_{0}(u))_{ss}}
e^{tY_{0}(u)-t(Y_{0}(u))_{ss}}v
=e^{t\lambda}e^{tY_{0}(u)-t(Y_{0}(u))_{ss}}v$. Since $V$ is spanned by such 
generalized eigenvector $v\in V$ for $Y_{0}(u)$, we have a well-defined 
operator $e^{tY_{0}(u)}$ on $V$ for every $t\in \C$.

\begin{prop}
For $t\in \C$, the operator $e^{tY_{0}(u)}$ is an automorphism of $V$.
\end{prop}
\pf
Since $Y_{0}(u)$ preserve the grading of $V$, 
$e^{tY_{0}(u)}$ also preserve the grading of $V$. 
{}From (\ref{comm-u0-ss}) and (\ref{comm-u0-nil}), we 
obtain 
$$e^{tY_{0}(u)}Y(v, x)=Y(e^{tY_{0}(u)}v, x)e^{tY_{0}(u)}.$$
The first half of the creation property says that $Y(u, x)\one
\in V[[x]]$. In particular, $Y_{0}(u)\one=0$. So $e^{tY_{0}(u)}\one=\one$.
Since $u\in V_{(1)}$ and $V_{(n)}=0$ for $n<0$, $L(n)u=0$ for $n\ge 2$.
We also have $L(1)u=0$. So $u$ is a lowest weight vector
of weight $1$ and $Y(u, x)$ is a primary field. 
Thus we have 
$$[L(-2), Y(u, x)]=x^{-1}\frac{d}{dx}Y(u, x)-z^{-2}Y(u, x).$$
Taking the coefficients of $x^{-1}$ of both sides of this formula, 
we obtain 
$$[L(-2), Y_{0}(u)]=0.$$
Thus $Y_{0}(u)\omega=L(-2)Y_{0}(u)\one=0$
which implies $e^{tY_{0}(u)}\omega=\omega$.
\epfv

By the skew-symmetry for $V$, we have 
$$Y(u, x)u=e^{xL(-1)}Y(u, -x)u.$$
Taking $\res_{x}$ of both sides, we have
$$Y_{0}(u)u=\sum_{m\in \N}\frac{(-1)^{-m-1}}{m!}(L(-1))^{m}Y_{m}(u)u.$$
Since $\wt Y_{m}(u)u=1-m-1+1=1-m<0$ and $V_{(n)}=0$ when $n<0$, 
$Y_{m}(u)u=0$ for $m>1$. Since $V_{(0)}=\C\mathbf{1}$,
$Y_{1}(u)u$ is proportional to $\mathbf{1}$ since 
$\wt Y_{1}(u)u=1-1-1+1=0$. So 
$(L(-1))^{m}Y_{m}(u)u=0$ for $m\ge 0$. Thus we obtain
$Y_{0}(u)u=-Y_{0}(u)u$ which implies $Y_{0}(u)u=0$.
By the component form of the
commutator formula (\ref{comm-u0}), we obtain
\begin{eqnarray*}
[Y_{0}(u), Y_{n}(u)]&=&Y_{n}((Y_{0}(u)u))\nn
&=&0
\end{eqnarray*}
for $n\in \Z$.
Then we can
define 
\begin{eqnarray*}
\Delta_{V}^{(u)}(x)&=&x^{Y_{0}(u)}e^{\int_{0}^{-x} Y^{\le -2}(u, y)}\nn
&=&x^{(Y_{0}(u))_{ss}}e^{(Y_{0}(u)-(Y_{0}(u))_{ss})\log x}
e^{\int_{0}^{-x} Y^{\le -2}(u, y)},
\end{eqnarray*}
where $(Y_{0}(u))_{ss}$ is the semisimple part of $Y_{0}(u)$ and, 
on an eigenvector of $(Y_{0}(u))_{ss}$ (that is, a generalized eigenvector
of $Y_{0}(u)$), 
$x^{(Y_{0}(u))_{ss}}$ is defined to be $x^{a}$ 
if the eigenvalue of the eigenvector is $a$.

\begin{thm}
For any $v\in V$, there exist $m_{1}, \dots, m_{l}\in \R$ such that 
$$\Delta_{V}^{(u)}(x)v\in x^{m_{1}}V[x^{-1}][\log x]+\cdots 
+x^{m_{l}}V[x^{-1}][\log x]$$
and the series $\Delta_{V}^{(u)}(x)$ satisfies 
\begin{equation}\label{Delta-thm-0-1}
\Delta_{V}^{(u)}(x)Y(v, x_{2})
=Y(\Delta_{V}^{(u)}(x+x_{2})v, x_{2})\Delta_{V}^{(u)}(x)
\end{equation}
and 
\begin{equation}\label{Delta-thm-0-2}
[L(-1), \Delta_{V}^{(u)}(x)]=-\frac{d}{dx}\Delta_{V}^{(u)}(x).
\end{equation}
\end{thm}
\pf
For any $v\in V$, 
\begin{eqnarray}\label{Delta-thm-1}
\lefteqn{e^{\int_{0}^{-x} Y^{\le -2}(u, y)v}v}\nn
&&=\sum_{k\in \N}\frac{1}{k!}
\left(\int_{0}^{-x} Y^{\le -2}(u, y)\right)^{k}v\nn
&&=\sum_{k\in \N}\frac{1}{k!}\left(\sum_{n\in \Z_{+}}\frac{(-1)^{n}}{-n}
Y_{n}(u) x^{-n}\right)^{k}v\nn
&&=\sum_{k\in \N}\frac{1}{k!}\sum_{n_{1}, \dots, n_{k}\in \Z_{+},\;
n_{1}+\cdots +n_{k}=l}
\frac{(-1)^{l}}{(-n_{1})\cdots (-n_{k})}
Y_{n_{1}}(u)\cdots Y_{n_{k}}(u)v x^{-l}.\nn
\end{eqnarray}
Since 
\begin{eqnarray*}
\wt Y_{n_{1}}(u)\cdots Y_{n_{k}}(u)
&=&(1-n_{1}-1)+\cdots +(1-n_{k}-1))\nn
&=&-n_{1}-\cdots -n_{k}\nn
&=&-l,
\end{eqnarray*}
we see that 
$Y_{n_{1}}(u)\cdots Y_{n_{k}}(u)w=0$ when 
$l$ is sufficiently large. Thus the right-hand side of 
(\ref{Delta-thm-1}) is an element of $V[x^{-1}]$ and
our first conclusion 
$$\Delta_{V}^{(u)}(x)v\in x^{m_{1}}V[x^{-1}][\log x]+\cdots 
+x^{m_{l}}V[x^{-1}][\log x]$$ 
follows immediately.

{}From (\ref{Delta-3}), we obtain
\begin{eqnarray}\label{Delta-thm-2}
\lefteqn{e^{\int_{0}^{-x} Y^{\le -2}(u, y)v}
Y(v, x_{2})e^{-\int_{0}^{-x} Y^{\le -2}(u, y)v}}\nn
&&=Y\left(e^{\int_{0}^{-x-x_{2}} Y^{\le -2}(u, y)v
+Y_{0}(u)\log\left(1+\frac{x_{2}}{x_{1}}\right)}v, x_{2}\right)\nn
&&=Y\left(e^{Y_{0}(u)\log\left(1+\frac{x_{2}}{x_{1}}\right)}
e^{\int_{0}^{-x-x_{2}} Y^{\le -2}(u, y)}v, x_{2}\right).
\end{eqnarray}

On the other hand, by (\ref{comm-u0-ss}) and 
(\ref{comm-u0-nil}), we have
\begin{eqnarray}\label{Delta-thm-3}
\lefteqn{x^{Y_{0}(u)}Y(v, x_{2})}\nn
&&=x^{(Y_{0}(u))_{ss}}e^{(Y_{0}(u)-(Y_{0}(u))_{ss})\log x}
Y(v, x_{2})\nn
&&=Y(x^{(Y_{0}(u))_{ss}}e^{(Y_{0}(u)-(Y_{0}(u))_{ss})\log x}v, x_{2})
x^{(Y_{0}(u))_{ss}}e^{(Y_{0}(u)-(Y_{0}(u))_{ss})\log x}\nn
&&=Y(x^{Y_{0}(u)}v, x_{2})x^{Y_{0}(u)}.
\end{eqnarray}
Using (\ref{Delta-thm-2}) and (\ref{Delta-thm-3}),
we obtain (\ref{Delta-thm-0-1}).

Finally, notice that both $[L(-1), \cdot]$ and $-\frac{d}{dx}$ are
derivatives
on the space $\C[[\int_{0}^{-x} Y^{\le -2}(u, y)]]$ of power series in 
$\int_{0}^{-x} Y^{\le -2}(u, y)$. So 
{}from (\ref{Delta-5}), we obtain
\begin{eqnarray}\label{Delta-thm-4}
\lefteqn{[L(-1), e^{\int_{0}^{-x} Y^{\le -2}(u, y)}]
+\frac{d}{dx}e^{\int_{0}^{-x} Y^{\le -2}(u, y)}}\nn
&&=\left[L(-1), \int_{0}^{-x} Y^{\le -2}(u, y)\right]
e^{\int_{0}^{-x} Y^{\le -2}(u, y)}\nn
&&\quad+\left(\frac{d}{dx}\int_{0}^{-x} Y^{\le -2}(u, y)\right)
e^{\int_{0}^{-x} Y^{\le -2}(u, y)}\nn
&&=-Y_{0}(u)x^{-1}e^{\int_{0}^{-x} Y^{\le -2}(u, y)}.
\end{eqnarray}
The equality (\ref{Delta-thm-4}) is equivalent to 
(\ref{Delta-thm-0-2}).
\epfv

Let $W=\coprod_{n, \alpha\in \C}
W_{\langle n \rangle}^{[\alpha]}$ 
be a strongly $\C$-graded 
generalized $V$-module (see \cite{HLZ1} and \cite{HLZ2})
with the action $e^{2\pi \sqrt{-1}\;(Y_{W})_{0}(u)}$
of the automorphism $e^{2\pi \sqrt{-1}\;Y_{0}(u)}$ of $V$
such that the $\C$-grading is given by the
generalized eigenspaces of
$(Y_{W})_{0}(u)$ and thus is compatible with the 
generalized eigenspaces of the action $e^{2\pi \sqrt{-1}\;(Y_{W})_{0}(u)}$
of $e^{2\pi \sqrt{-1}\;Y_{0}(u)}$. 
For example, when $W$ is a $V$-module, $W$ has a
natural $\C$-grading given by the generalized eigenspaces
for $(Y_{W})_{0}(u)$ and together with this grading,
$W$ is such a generalized $V$-module.
We define a linear map
\begin{eqnarray*}
Y_{W}^{(u)}: V\otimes W&\to&  W\{x\}[\log x]\nn
v\otimes w&\mapsto& Y_{W}^{(u)}(v, x)w
\end{eqnarray*}
by 
$$Y_{W}^{(u)}(v, x)w=Y_{W}(\Delta_{V}^{(u)}(x)v, x).$$

First, we have:

\begin{lemma}\label{new-l-0}
The vertex operator $Y_{W}^{(u)}(\omega, x)$ is in $W[[x, x^{-1}]]$
and if we write 
$$Y_{W}^{(u)}(\omega, x)=\sum_{n\in \Z}L_{W}^{(u)}(n)x^{-n-2},$$
then 
$$L_{W}^{(u)}(0)=L_{W}(0)+(Y_{W})_{0}(u)-\frac{1}{2}\mu,$$
where 
$\mu\in \C$ is given by 
$$(Y_{W})_{1}(u)u=\mu \one$$
(note that since
$(Y_{W})_{1}(u)u\in V_{(0)}=\C\one$,
it must be proportional to 
$\one$). In particular, for $n\in \C$
and $\alpha$ an eigenvalue of $(Y_{W})_{0}(u)$,
$W_{\langle n-\alpha+\frac{1}{2}\mu \rangle}^{[\alpha]}$ 
is the generalized eigenspace
for both $L_{W}^{(u)}(0)$ and $(Y_{W})_{0}(u)$ 
with the eigenvalues
$n$ and $\alpha$, respectively.
\end{lemma}
\pf
By (\ref{Delta-1}), we have
$${[Y(\omega, x_{1}), Y(u, x_{2})]}=
\res_{x_{0}}x^{-1}_{2}\delta\left(\frac{x_{1}-x_{0}}{x_{2}}\right)
Y(Y(\omega, x_{0})u, x_{2}).$$
Taking the coefficients of $x_{1}^{0}$ in both sides and noticing that from 
our assumption, $L(n)u=0$ for $n>0$, we obtain
\begin{eqnarray*}
[L(-2), Y(u, x_{2})]&=&
-x_{2}^{-1}Y(L(-1)u, x_{2})+x_{2}^{-2}Y(L(0)u, x_{2})\nn
&=&-x_{2}^{-1}\frac{d}{dx_{2}}Y(u, x_{2})+x_{2}^{-2}Y(u, x_{2}).
\end{eqnarray*}
Taking the coefficients of $x_{2}^{-m-1}$ in both sides,
we obtain
$$[L(-2), Y_{m}(u)]=mY_{m-2}(u)$$
for $m\in \Z$. 
Then we have 
\begin{eqnarray*}
Y_{m}(u)\omega&=&Y_{m}(u)L(-2)\one\nn
&=&L(-2)Y_{m}(u)\one-mY_{m-2}(u)\one
\end{eqnarray*}
for $m\in \Z$. Since $Y_{m}(u)\one=0$ for $m\ge 0$ and $Y_{-1}(u)\one=u$,
we obtain
$$Y_{m}(u)\omega=0$$
for $m\in\N$ but $m\ne 1$ and
$$Y_{1}(u)\omega=-u.$$
In particular, $Y_{0}(u)\omega=(Y_{0}(u))_{ss}\omega$.
We also know that 
$Y_{m}(u)u=0$ for $m>1$ and $Y_{1}(u)u=\mu \one$.
So 
$$e^{\int_{0}^{-x} Y^{\le -2}(u, y)}\omega
=\omega+ux^{-1}-\frac{1}{2}\mu\one x^{-2}$$
and 
\begin{eqnarray*}
\Delta_{V}^{(u)}(x)\omega&=&x^{(Y_{0}(u))_{ss}}
(\omega+ux^{-1}-\frac{1}{2}\mu\one x^{-2})\nn
&=&\omega+ux^{-1}-\frac{1}{2}\mu\one x^{-2}.
\end{eqnarray*}
Thus
$$Y_{W}^{(u)}(\omega, x)=Y_{W}(\omega, x)
+Y(u, x)x^{-1}-\frac{1}{2}\mu x^{-2}\in W[[x, x^{-1}]]$$
and
$L_{W}^{(u)}(0)=L_{W}(0)+(Y_{W})_{0}(u)-\frac{1}{2}\mu.$

The second conclusion follows immediately.
\epfv

We have the following consequence:

\begin{cor}\label{double-grading}
The space $W$ is also doubly graded by the generalized 
eigenspaces for $L_{W}^{u}(0)$ and $(Y_{W})_{0}(u)$, that is, 
$W$ has a $\C\times \C$-grading 
$W=\coprod_{n, \alpha\in \C}W_{[n]}^{[\alpha]}$
where for $n, \alpha\in \C$,
$$W_{[n]}^{[\alpha]}
=W_{\langle n-\alpha+\frac{1}{2}\mu \rangle}^{[\alpha]},$$
that is, 
$w\in W_{[n]}^{[\alpha]}$ for $n, \alpha\in \C$
if and only if $w$ is a generalized eigenvector 
for $L_{W}^{u}(0)$ with the eigenvalue $n$ and a generalized eigenvector 
for $(Y_{W})_{0}(u)$ with the eigenvalue $\alpha$.
\end{cor}
\pf
By Lemma \ref{new-l-0},
every element of $W$ is a finite sum of 
vectors which are generalized eigenvectors for both
$L_{W}^{(u)}(0)$ and  $(Y_{W})_{0}(u)$.
\epfv

The following theorem is the main result of the present 
paper:

\begin{thm}\label{main}
The pair $(W, Y_{W}^{(u)})$, where $W$ is equipped with the 
$\C\times \C$-grading given in Corollary \ref{double-grading}, 
is a strongly $\C$-graded generalized
$e^{2\pi \sqrt{-1}\; Y_{0}(u)}$-twisted $V$-module.
\end{thm}
\pf
We first prove that 
$W$ satisfies the grading-restriction 
condition. By Lemma \ref{new-l-0},
$L_{W}(0)=L_{W}^{(u)}(0)-(Y_{W})_{0}(u)+\frac{1}{2}\mu$. 
Since the generalized $V$-module $W$ satisfies the
grading-restriction condition, 
$\dim W_{\langle n \rangle}^{[\alpha]}<\infty$ for $n, \alpha\in \C$
and for any fixed $n, \alpha\in \C$,
$W_{\langle n+l \rangle}^{[\alpha]}=0$ when $l$ is a 
sufficiently negative integer. 
So for $n, \alpha\in \C$,
$W_{[n]}^{[\alpha]}
= W_{\langle n-\alpha+\frac{1}{2}\mu \rangle}^{[\alpha]}$
are finite dimensional and for any fixed $n, \alpha\in \C$,
$W_{[n+l]}^{[\alpha]}= 
W_{\langle n+l-\alpha+\frac{1}{2}\mu \rangle}^{[\alpha]}=0$
when $l$ is a 
sufficiently negative integer.

{}From the first part of the creation property,
$Y(u, x)\one\in V[[x]]$ for vertex operator algebra, we have 
$Y_{n}(u)\one=0$ 
for $n\ge 0$. In particular, 
$(Y_{0}(u))_{ss}\one=(Y_{0}(u)-(Y_{0}(u))_{ss})\one
=0$. We obtain
\begin{eqnarray*}
\Delta_{V}^{(u)}(x)\one
&=&x^{(Y_{0}(u))_{ss}}e^{(Y_{0}(u)-(Y_{0}(u))_{ss})\log x}
e^{\int_{0}^{-x} Y^{\le -2}(u, y)}\one\nn
&=&\one.
\end{eqnarray*}
For $w\in W$, 
\begin{eqnarray*}
Y_{W}^{(u)}(\one, x)w&=&Y_{W}(\Delta_{V}^{(u)}(x)\one, x)w\nn
&=&Y_{W}(\one, x)w\nn
&=&w,
\end{eqnarray*}
proving the identity property.

Let $v\in V$ be a generalized eigenvector for $Y_{0}(u)$ with 
eigenvalue $\lambda$. Then for $v\in V$ and $z\in \C^{\times}$,
\begin{eqnarray*}
\lefteqn{\Delta_{V}^{(u)}(x)e^{2\pi i
Y_{0}(u)}v\lbar_{x^{n}=e^{nl_{p}(z)},\;\log x=l_{p}(z)}}\nn
&&=x^{\lambda}e^{(Y_{0}(u)-(Y_{0}(u))_{ss})\log x}
e^{2\pi i\lambda}e^{2\pi i (Y_{0}(u)-(Y_{0}(u))_{ss})}
v\lbar_{x^{n}=e^{nl_{p}(z)},\;\log x=l_{p}(z)}\nn
&&=x^{\lambda}e^{(Y_{0}(u)-(Y_{0}(u))_{ss})\log x}
v\lbar_{x^{n}=e^{nl_{p+1}(z)},\;\log x=l_{p+1}(z)}\nn
&&=\Delta_{V}^{(u)}(x)v\lbar_{x^{n}=e^{nl_{p+1}(z)},\;\log x=l_{p+1}(z)}.
\end{eqnarray*}
Thus
\begin{eqnarray*}
\lefteqn{(Y_{W}^{(u)})^{e^{2\pi iY_{0}(u)}; p} (e^{2\pi iY_{0}(u)}v, z)w}\nn
&&=Y_{W}(\Delta_{V}^{(u)}(x)e^{2\pi i
Y_{0}(u)}v, x)w\lbar_{x^{n}=e^{nl_{p}(z)},\;\log x=l_{p}(z)}\nn
&&=Y_{W}(\Delta_{V}^{(u)}(x)v, x)w\lbar_{x^{n}=e^{nl_{p+1}(z)},\;
\log x=l_{p+1}(z)}\nn
&&= (Y_{W}^{(u)})^{e^{2\pi iY_{0}(u)}; p+1}(v, z)w
\end{eqnarray*}
for $v\in V$ and $w\in W$, that is, the equivariance 
property holds.

By Corollary \ref{double-grading},
the $L^{(u)}(0)$-grading condition and the grading-compatibility
condition are satisfied. 

By the $L(-1)$-derivative property for the vertex operator 
map $Y_{W}$ and (\ref{Delta-thm-0-2}), we have 
\begin{eqnarray*}
\frac{d}{dx}Y_{W}^{(u)}(v, x)&=&
\frac{d}{dx}Y_{W}(\Delta_{V}^{(u)}(x)v, x)\nn
&=&Y_{W}(L(-1)\Delta_{V}^{(u)}(x)v, x)
+Y_{W}\left(\left(\frac{d}{dx}\Delta_{V}^{(u)}(x)v\right), x\right)\nn
&=&Y_{W}(L(-1)\Delta_{V}^{(u)}(x)v, x)
+Y_{W}(-[L(-1), \Delta_{V}^{(u)}(x)]v), x)\nn
&=&Y_{W}(\Delta_{V}^{(u)}(x)L(-1)v, x)\nn
&=&Y_{W}^{(u)}(L(-1)v, x),
\end{eqnarray*}
for $v\in V$, proving the $L(-1)$-derivative property. 

Finally we prove the duality property. 
By Corollary \ref{double-grading}, 
$W_{[n]}^{[\alpha]}=W_{\langle n-\alpha+\frac{1}{2}\mu \rangle}^{[\alpha]}$.
So the graded dual of $W$ with 
respect to the grading 
$W=\coprod_{n, \alpha\in \C}W_{\langle n \rangle}^{[\alpha]}$
and the graded dual with respect the new $\C\times \C$-grading 
$W=\coprod_{n, \alpha\in \C}W_{[n]}^{[\alpha]}$
in Corollary \ref{double-grading}
are equal as vector spaces, though their gradings are 
different. We shall use $W'$ to denote the underlying 
vector space of these two graded duals. 
It will be clear from the context  which grading we will
be using.

For $v_{1}, v_{2}\in V$, $w\in W$ and $w'\in W'$, we know that
there exist $m_{1}, \dots, m_{r}$, $n_{1}, \dots, n_{s}\in \R$
such that 
$$\Delta_{V}^{(u)}(x_{1})e^{2\pi i
Y_{0}(u)}v_{1}\in x_{1}^{m_{1}}V[x_{1}^{-1}][\log x_{1}]+\cdots 
+x_{1}^{m_{r}}V[x_{1}^{-1}][\log x_{1}]$$ 
and 
$$\Delta_{V}^{(u)}(x_{2})e^{2\pi i
Y_{0}(u)}v_{2}\in x_{2}^{n_{1}}V[x_{2}^{-1}][\log x_{2}]+\cdots 
+x_{2}^{n_{s}}V[x_{2}^{-1}][\log x_{2}].$$
Thus, using
the rationality, commutativity and associativity properties for 
the $V$-module $W$ and the fact that in the region 
$|z_{2}|>|z_{1}-z_{2}|>0$
\begin{eqnarray*}
\lefteqn{\Delta_{V}^{(u)}(x_{2}+x_{0})v_{1}
\lbar_{x_{0}^{n}=(z_{1}-z_{2})^{n},\;x_{2}^{n}=e^{nl_{p}(z_{2})},\;
\log x_{2}=l_{p}(z_{2})}}\nn
&&=\Delta_{V}^{(u)}(x_{1})v_{1}\lbar_{x_{1}^{n}=e^{nl_{p}(z_{1})},\;
\log x_{1}=l_{p}(z_{1})},
\end{eqnarray*}
we see that
there exists $a_{ijkl}\in \C$ for $i, j, k, l=0, \dots, N$,
$\alpha_{i}, \beta_{j}\in \C$ for $i, j=0, \dots, N$
and $t\in \N$ such that
the series 
\begin{eqnarray*}
\lefteqn{\langle w', (Y_{W}^{(u)})^{e^{2\pi iY_{0}(u)}; p}(v_{1}, z_{1})
(Y_{W}^{(u)})^{e^{2\pi iY_{0}(u)}; p}(v_{2}, z_{2})w\rangle}\nn
&&=\langle w', Y_{W}(\Delta_{V}^{(u)}(x_{1})
v_{1}, x_{1})\cdot\nn
&&\quad\quad\quad
\cdot Y_{W}(\Delta_{V}^{(u)}(x_{2})v_{2}, x_{2})w
\rangle\lbar_{x_{1}^{n}=e^{nl_{p}(z_{1})},\;
\log x_{1}=l_{p}(z_{1}),\;x_{2}^{n}=e^{nl_{p}(z_{2})},\;
\log x_{2}=l_{p}(z_{2})},\\
\lefteqn{\langle w', (Y_{W}^{(u)})^{e^{2\pi iY_{0}(u)}; p}(v_{2}, z_{2})
(Y_{W}^{(u)})^{e^{2\pi iY_{0}(u)}; p}(v_{1}, z_{1})w\rangle}\nn
&&=\langle w', Y_{W}(\Delta_{V}^{(u)}(x_{2})
v_{2}, x_{2})\cdot\nn
&&\quad\quad\quad
\cdot Y_{W}(\Delta_{V}^{(u)}(x_{1})
v_{1}, x_{1})w\rangle\lbar_{x_{1}^{n}=e^{nl_{p}(z_{1})},\;
\log x_{1}=l_{p}(z_{1}),\;x_{2}^{n}=e^{nl_{p}(z_{2})},\;
\log x_{2}=l_{p}(z_{2})},\\
&&\langle w', Y_{W}(Y(\Delta_{V}^{(u)}(x_{2}+x_{0})v_{1}, x_{0})
\cdot\nn
&&\quad\quad\quad\quad\quad\quad\quad\quad\quad
\cdot \Delta_{V}^{(u)}(x_{2})v_{2}, x_{2})w\rangle
\lbar_{x_{0}^{n}=(z_{1}-z_{2})^{n},\;x_{2}^{n}=e^{nl_{p}(z_{2})},\;
\log x_{2}=l_{p}(z_{2})}
\end{eqnarray*}
are absolutely convergent in the regions $|z_{1}|>|z_{2}|>0$,
$|z_{1}|>|z_{2}|>0$, $|z_{2}|>|z_{1}-z_{2}|>0$, respectively,
to the branch 
$$\sum_{i, j, k, l=0}^{N}
a_{ijkl}e^{m_{i}l_{p}(z_{1})}e^{n_{j}l_{p(z_{2})}}
l_{p}(z_{1})^{k}l_{p}(z_{2})^{l}(z_{1}-z_{2})^{-t}$$
of the multivalued analytic function 
$$\sum_{i, j, k, l=0}^{N}a_{ijkl}
z_{1}^{m_{i}}z_{2}^{n_{j}}(\log z_{1})^{k}
(\log z_{2})^{l}(z_{1}-z_{2})^{-t}.$$
By (\ref{Delta-thm-0-1}), 
\begin{eqnarray*}
\lefteqn{\langle w', Y_{W}(Y(\Delta_{V}^{(u)}(x_{2}+x_{0})v_{1}, x_{0})
\cdot}\nn
&&\quad\quad\quad\quad\quad\quad\quad\quad
\cdot \Delta_{V}^{(u)}(x_{2})v_{2}, x_{2})w\rangle
\lbar_{x_{0}^{n}=(z_{1}-z_{2})^{n},\;x_{2}^{n}=e^{nl_{p}(z_{2})},\;
\log x_{2}=l_{p}(z_{2})}\nn
&&=\langle w', Y_{W}(\Delta_{V}^{(u)}(x_{2})Y(v_{1}, x_{0})
v_{2}, x_{2})w\rangle
\lbar_{x_{0}^{n}=(z_{1}-z_{2})^{n},\;x_{2}^{n}=e^{nl_{p}(z_{2})},\;
\log x_{2}=l_{p}(z_{2})}\nn
&&=\langle w', (Y_{W}^{(u)})^{e^{2\pi iY_{0}(u)}; p}
(Y(v_{1}, z_{1}-z_{2})
v_{2}, z_{2})w\rangle.
\end{eqnarray*}
Thus the duality property is proved.
\epfv

\begin{rema}
{\rm In the case that $\res_{x} Y(u, x)$
acts on $V$ semisimply and has eigenvalues belonging to $\frac{1}{k}\Z$,
the theorem above reduces to the construction by Li in
\cite{Li}.}
\end{rema}

We have the following immediate consequence:

\begin{cor}
The correspondence given by $(W, Y_{W})\to (W, Y_{W}^{(u)})$ 
is a functor from the category of strongly $\C$-graded
generalized $V$-modules to the category of strongly $\C$-graded
generalized $e^{2\pi \sqrt{-1} \; Y_{0}(u)}$-twisted $V$-modules. \epf
\end{cor}

We now apply the theorem above to give some explicit examples. 
Let $p$ and $q$ be a pair of coprime positive integers,
$L$ the one-dimensional positive-definite even 
lattice generated by $\gamma$ with the 
bilinear form $\langle\cdot, \cdot\rangle$ given by 
$\langle \gamma, \gamma\rangle=2pq$, and 
$\widetilde{L}$ its dual lattice. 
Let $V_{L}$ be the vertex algebra associated to $L$ and 
$V_{\widetilde{L}}$ the generalized vertex algebra associated
to $\widetilde{L}$ (see \cite{DL}).
The element
$$\omega=\frac{1}{pq}\gamma(-1)^{2}\one+\frac{p-q}{2pq}\gamma(-2)\one
\in V_{L}\subset V_{\widetilde{L}}$$
is a conformal element. If instead of the usual
conformal element $\frac{1}{pq}\gamma(-1)^{2}\one$,
we take $\omega$ to be the 
conformal element for $V_{L}$ and $V_{\widetilde{L}}$, we obtain 
a vertex operator algebra (still denoted $V_{L}$) and a generalized
vertex operator algebra (still denoted $V_{\widetilde{L}}$), respectively.
Note that in the grading given by the conformal element 
$\omega$, the weights of $e^{\gamma/q}$ and $e^{-\gamma/p}$ are $1$. 
Consider the $V_{L}$-modules $V_{L-\frac{\gamma}{p}}$
and $V_{L+\frac{\gamma}{q}}$ which are both graded 
by $\Z$. Then we have vertex-operator-algebraic extensions 
$\mathcal{V}_{0}=V_{L}\oplus V_{L-\frac{\gamma}{p}}$,
$\mathcal{V}_{0}^{o}=V_{L}\oplus V_{L+\frac{\gamma}{q}}$
and $\mathcal{V}(p, q)=V_{L}\oplus V_{L-\frac{\gamma}{p}}
\oplus V_{L+\frac{\gamma}{q}}$ of $V_{L}$
for which the vertex operators are given by
\begin{eqnarray*}
Y_{\mathcal{V}_{0}}(u_{1}+v_{1}, x)(u_{2}+v_{2})
&=&Y_{V_{\widetilde{L}}}(u_{1}+v_{1}, x)u_{2}
+Y_{V_{\widetilde{L}}}(u_{1}, x)v_{2},\\
Y_{\mathcal{V}_{0}^{o}}(u_{1}+w_{1}, x)(u_{2}+w_{2})
&=&Y_{V_{\widetilde{L}}}(u_{1}+w_{1}, x)u_{2}
+Y_{V_{\widetilde{L}}}(u_{1}, x)w_{2},
\end{eqnarray*}
\vspace{-3em}
\begin{eqnarray*}
\lefteqn{Y_{\mathcal{V}(p, q)}(u_{1}+v_{1}+w_{1}, x)(u_{2}+v_{2}+w_{2})}\nn
&&\quad\quad\quad\quad\quad\quad\quad\quad\quad\quad
=Y_{V_{\widetilde{L}}}(u_{1}+v_{1}+w_{1}, x)u_{2}
+Y_{V_{\widetilde{L}}}(u_{1}, x)(v_{2}+w_{2})
\end{eqnarray*}
for $u_{1}, u_{2}\in V_{L}$ and $v_{1}, v_{2}\in 
V_{L-\frac{\gamma}{p}}$ and  $w_{1}, w_{2}\in 
V_{L+\frac{\gamma}{q}}$.
Note that $\mathcal{V}_{0}$ and $\mathcal{V}_{0}^{o}$
are vertex operator subalgebras of $\mathcal{V}(p, q)$ and that
$e^{-\gamma/p}\in \mathcal{V}_{0}$ and $e^{\gamma/q}\in \mathcal{V}_{0}^{o}$.
Also note that $(\mathcal{V}_{0})_{(0)}=(\mathcal{V}_{0}^{o})_{(0)}
=(\mathcal{V}(p, q))_{(0)}=\C\one$ and 
$(\mathcal{V}_{0})_{(n)}=(\mathcal{V}_{0}^{o})_{(n)}
=(\mathcal{V}(p, q))_{(n)}=0$ for $n<0$.

Let 
$Q= \res_x Y_{\mathcal{V}(p, q)}(e^{\gamma/q}, x)$ and $\tilde{Q} = \res_{x}
Y_{\mathcal{V}(p, q)}(e^{-\gamma/p}, x)$. These
operators are called {\it screening operators}.
Then $e^{2\pi \sqrt{-1}\;\tilde{Q}}$ and $e^{2\pi \sqrt{-1}\;Q}$ 
are automorphisms of $\mathcal{V}_{0}$ and $\mathcal{V}_{0}^{o}$, respectively,
and are both automorphisms of $\mathcal{V}(p, q)$. 
In fact, the triplet vertex operator algebra $\mathcal{W}(p, q)$
(see \cite{FGST} and \cite{AM}) is a vertex operator subalgebra
of the fixed-point subalgebra
of $\mathcal{V}(p, q)$ 
under the group generated by $e^{2\pi \sqrt{-1}\;Q}$ and 
$e^{2\pi \sqrt{-1}\;\tilde{Q}}$.

Let $W=\coprod_{n, \alpha\in \C}
W_{\langle n \rangle}^{[\alpha]}$ 
be a strongly $\C$-graded 
generalized $\mathcal{V}_{0}$-, $\mathcal{V}_{0}^{o}$-
or $\mathcal{V}(p, q)$-module
with an action $e^{2\pi \sqrt{-1}\;(Y_{W})_{0}(e^{\gamma/q})}$ of  
$e^{2\pi \sqrt{-1}\;\tilde{Q}}$, 
$e^{2\pi \sqrt{-1}\;(Y_{W})_{0}(e^{-\gamma/p})}$ of
$e^{2\pi \sqrt{-1}\;Q}$ or either of them, respectively,
such that the $\C$-grading is given by
the generalized eigenspaces of $(Y_{W})_{0}(e^{\gamma/q})$
or $(Y_{W})_{0}(e^{-\gamma/p})$. For example, 
we can take $W$ to be $\mathcal{V}_{0}$, $\mathcal{V}_{0}^{o}$
or $\mathcal{V}(p, q)$ themselves.
In \cite{AM} (Theorem 9.1), Adamovi\'{c} and Milas in \cite{AM} proved 
that $Y_{W}^{(e^{-\gamma/p})}$ and $Y_{W}^{(e^{\gamma/q})}$
are intertwining operators of suitable types. Applying Theorem
\ref{main}, we obtain immediately:

\begin{thm}
For a strongly $\C$-graded 
generalized $\mathcal{V}_{0}$-module (or $\mathcal{V}_{0}^{o}$-
or $\mathcal{V}(p, q)$-module) $(W, Y_{W})$,
the pair $(W, Y_{W}^{(e^{\gamma/q})})$ (or the pair $(W, Y_{W}^{(e^{-\gamma/p})})$
or the pairs $(W, Y_{W}^{(e^{\gamma/q})})$ and $(W, Y_{W}^{(e^{-\gamma/p})})$)
is a strongly $\C$-graded 
generalized $e^{2\pi \sqrt{-1}\;\tilde{Q}}$-twisted 
$\mathcal{V}_{0}$-modules (or is a strongly $\C$-graded 
generalized $e^{2\pi \sqrt{-1}\; Q}$-twisted 
$\mathcal{V}_{0}^{o}$-module or are strongly $\C$-graded 
generalized $e^{2\pi \sqrt{-1}\;\tilde{Q}}$- and $e^{2\pi \sqrt{-1}\; Q}$-twisted
$\mathcal{V}(p, q)$-modules, respectively). \epf
\end{thm}

\noindent {\small \sc Department of Mathematics, Rutgers University,
110 Frelinghuysen Rd., Piscataway, NJ 08854-8019}

\noindent {\em E-mail address}: yzhuang@math.rutgers.edu

\end{document}